\documentclass[11pt]{article}

\usepackage{amsfonts}
\usepackage{epsfig}
\usepackage{amsmath}
\numberwithin{equation}{section}
\usepackage{hyperref}
\usepackage{amssymb}
\usepackage{graphics}
\usepackage{verbatim}
\usepackage{psfrag}
\usepackage{subfigure}
\usepackage{color}
\usepackage{amsthm}
\usepackage[all]{xy}
\usepackage{makeidx}
\usepackage{enumerate}
\usepackage[paperheight=235mm,paperwidth=155mm,textwidth=125mm,textheight=199mm]{geometry}

\let\oldbibliography\thebibliography
\renewcommand{\thebibliography}[1]{%
  \oldbibliography{#1}%
  \setlength{\itemsep}{-1.2mm}%
}

\theoremstyle{plain}
\newtheorem{thm}{Theorem}[section]

\theoremstyle{definition}
\newtheorem{defn}[thm]{Definition}

\newtheorem{ex}[thm]{Example}

\newtheoremstyle{myremark}
  {3pt}
  {3pt}
  {\small \rmfamily}
  {5pt}
  {\rmfamily}
  {:}
  {.5em}
  {}

\theoremstyle{myremark}
\newtheorem*{remark}{\textit{Remark}}

\def\R{\mathbb{R}}

\def\N{\mathbb{N}}
\def\Z{\mathbb{Z}}
\def\P{\mathbb{P}}
\def\E{\mathbb{E}}
\def\T{\mathbb{T}}

\def\cA{\mathcal{A}}
\def\cB{\mathcal{B}}

\def\cD{\mathcal{D}}

\def\cF{\mathcal{F}}
\def\cG{\mathcal{G}}
\def\cH{\mathcal{H}}
\def\cI{\mathcal{I}}

\def\cK{\mathcal{K}}

\def\cM{\mathcal{M}}
\def\cN{\mathcal{N}}
\def\cO{\mathcal{O}}
\def\cP{\mathcal{P}}

\def\cR{\mathcal{R}}
\def\cS{\mathcal{S}}
\def\cT{\mathcal{T}}
\def\cU{\mathcal{U}}
\def\cV{\mathcal{V}}
\def\cW{\mathcal{W}}

\def\txtc{{\textnormal{c}}}
\def\txtd{{\textnormal{d}}}
\def\txte{{\textnormal{e}}}

\def\txts{{\textnormal{s}}}

\def\txtD{{\textnormal{D}}}

\def\Id{{\textnormal{Id}}}

\def\ra{\rightarrow}
\def\I{\infty}
\def\fs{\mathfrak{s}}

\newcommand{\be}{\begin{equation}}
\newcommand{\ee}{\end{equation}}
\newcommand{\benn}{\begin{equation*}}
\newcommand{\eenn}{\end{equation*}}
\newcommand{\bea}{\begin{eqnarray}}
\newcommand{\eea}{\end{eqnarray}}
\newcommand{\beann}{\begin{eqnarray*}}
\newcommand{\eeann}{\end{eqnarray*}}

\newcommand{\myendex}{$\blacklozenge$\end{ex}}
\newcommand{\myendexerc}{$\lozenge$\end{exerc}}
\newcommand{\myendpexerc}{$\lozenge$\end{pexerc}}

\def\XXint#1#2#3{{\setbox0=\hbox{$#1{#2#3}{\int}$}
\vcenter{\hbox{$#2#3$}}\kern-.5\wd0}}

\makeindex

\begin{document}

\author{
$\quad$\\
$\quad$\\
\Large{Christian Kuehn \& Alexandra Neam\c tu}\\
$\quad$\\
$\quad$\\
$\quad$\\
$\quad$\\
Faculty of Mathematics,\\ 
Technical University of Munich, \\
85748 Garching~bei~M\"unchen, Germany}

\title{Dynamics of Stochastic \\ Reaction-Diffusion Equations}

\maketitle

\section{Introduction}
Stochastic partial differential equations (SPDEs) represent a very active research field with numerous recent developments and breakthrough results~\cite{Hairer,Hairer2}.
There are several well-established approaches and methods used to construct solutions for SPDEs, which is always a challenge due to the irregularity of the noise terms that perturb the equation. In applications, such noise terms can quantify the lack of knowledge of certain parameters, finite-size effects, and/or fluctuations occurring due to external perturbations.
Since SPDEs have become a key modelling tool in applications, there has been a growing
interest in studying their dynamical phenomena~\cite{crauel, CerraiRoeckner, DLS, HairerE, BerglundGentzM, MuellerMQ}.\\
The main goal of this work is to provide a survey on different approaches to solution theory and dynamical properties for SPDEs, which is accessible for a wide community interested in modern methods in stochastic analysis, dynamics and applications. \\
This work is structured as follows: Section~\ref{chapter2} is devoted to the solution theory of SPDEs driven by general additive and multiplicative noise.~The aim is to illustrate the classical It\^o-theory, the random field approach and also recent developments in the context of rough paths and regularity structures. The first three subsections~\ref{two:1}-\ref{variational} provide a concise survey of several broadly used classical solution concepts for semi-linear SPDEs such as mild, weak, strong, martingale and variational solutions~\cite{DaPratoZabczyk,PrevotRoeckner,KrylovR}. In Section~\ref{quasilinear} we investigate more complicated quasilinear SPDEs and their pathwise mild solutions~\cite{KuehnNeamtu}. Section~\ref{random:field} presents an alternative to the It\^o calculus via the random field approach~\cite{Dalang, DalangLluis} based on Walsh integration theory~\cite{Walsh}.
In contrast to the first five subsections of this chapter, Section~\ref{rough:paths} starts by presenting one possible pathwise approach to the construction of stochastic integrals and solutions of SPDEs. Pathwise techniques are not restricted to semi-martingales as the It\^o calculus or Walsh theory and apply, for instance, to certain types of fractional Brownian motion~\cite{Gubinelli,GarridoLuSchmalfuss,HesseNeamtu1}. Finally, we conclude the part on solution theory in Section~\ref{chapter2} with singular SPDEs. Due to the high irregularity of the noise, such SPDEs become ill-posed, as also justified in Section~\ref{random:field}. Therefore, new tools such as regularity structures~\cite{Hairer} are required in order to give an appropriate meaning of the solution.\\
In Section~\ref{dynamics} the emphasis is on the long-time qualitative and quantitative dynamical behavior of SPDEs. Section~\ref{rds} collects concepts from the theory of random dynamical systems~\cite{Arnold} which will be employed in Sections~\ref{im} and~\ref{ra} to construct random invariant manifolds~\cite{DuanLuSchmalfuss, Bloemker, GS, KN, WanngDuan} and attractors~\cite{crauel,CrauelFlandoli,Debussche,schmalfuss, Gess1, Gess2, GessLiuRoeckner} for SPDEs. We provide examples for SPDEs which fit into the theory of random dynamical systems and discuss the difficulties that occur in this sense for general SPDEs driven by multiplicative noise. To overcome such obstacles it is helpful to rely on a pathwise approach. Based on this we provide a center manifold theory for PDEs with irregular forcing in Section~\ref{im}. Section~\ref{stability} deals with various stability concepts for SPDEs such as stability in probability, almost sure exponential stability and moment exponential stability~\cite{Mao, DeyaStability, CaraballoLiu, GS}. We point out the main techniques required in this setting and illustrate them on several SPDE examples. We do not focus only on asymptotic stability, but also discuss metastability, i.e.~dynamics which looks stable on very long time scales~\cite{BerglundGentzM, Barret, BerglundAC}. Here we explain, how large deviation principles~\cite{FreidlinWentzell, CerraiRoeckner, HSS} can be helpful in order to investigate such a behavior. We conclude Section~\ref{dynamics} mentioning further dynamical phenomena for SPDEs, which were not discussed in this work due to space limitations, and we provide a brief outlook.


\section{Solution Theory for SPDEs}
\label{chapter2}

We start by fixing the main notation used throughout this chapter. Let $(\cH,\|\cdot\|_{\cH})$ be an infinite-dimensional separable Hilbert space and let $\mathcal{L}(\cH)$ denote the space of linear operators on $\cH$ and $\mathcal{L}_{2}(\cH)$ the space of Hilbert-Schmidt operators. The scalar product on a Hilbert space will be denoted by $[\cdot,\cdot]$ and $\left<\cdot,\cdot\right>$ stands for the duality pairing. The dual space of a Banach space $\cV$ will be denoted by $\cV^{*}$.~For a linear operator $A$, we write $A^{*}$ for its adjoint and $\textnormal{Dom}(A)$ for its domain.~Furthermore, $(\Omega,\mathcal{F},\mathbb{P})$ denotes a complete probability space, $(\mathcal{F}_{t})_{t\geq 0}$ is a filtration, $\mathcal{B}(\cdot)$ signifies the Borel-sigma algebra on a certain space, $\otimes$ indicates the product $\sigma$-algebra as well as product measures, and $\mathbb{E}$ stands for the mathematical expectation. We use the abbreviations $\mathbb{P}$-a.s.~(almost surely) and a.e.~(almost everywhere). For $d\geq 1$ and a domain $D\subseteq\mathbb{R}^d$ and $1<p<\infty$, $L^{p}(D,\cH)$ denotes the space of $p$-integrable (over $D$) $\cH$-valued functions, whereas $L^{p}(\Omega,\mathcal{F},\mathbb{P};\cH)$ refers to the space of $\cH$-valued random variables defined on the probability space $(\Omega,\mathcal{F},\mathbb{P})$ having finite moments of order $p$. For $k\in\mathbb{R}$ and $1<p<\infty$, $W^{k,p}(D,\cH)$ denote the usual Sobolev spaces (i.e.~contain functions that belong to $L^{p}(D,\cH)$ and their weak derivatives up to order $k$ are again in $L^{p}(D,\cH)$). Whenever we incorporate Dirichlet / Neumann boundary conditions, we write $W^{k,p}_{D}$ / $W^{k,p}_{N}$. For zero Dirichlet boundary conditions we write for simplicity $W^{k,p}_0$. Finally we fix a time horizon $T>0$ and set $\Delta_{T}:=\{(s,t)\in[0,T]\times[0,T]: t\geq s\}$. \\
 In order to deal with stochastic evolution equations in infinite-dimensional spaces we firstly introduce the following Hilbert space-valued processes. Further details can be looked up in~\cite[Chapter~4]{DaPratoZabczyk},~\cite[Chapter~6]{isem} and~\cite{PrevotRoeckner, KrylovR}.

\begin{defn}
	An {\bf $\cH$-isonormal process} on $\Omega$ is a mapping $\cW:\cH\to L^{2}(\Omega)$ such that the following two conditions hold:
	\begin{itemize}
		\item for every $h\in\cH$ the random variable $\cW h$ is centered Gaussian, more precisely $\cW h\sim N(0,\|h\|^2_{\cH})$;
		\item $\mathbb{E}(\cW h_1\cdot \cW h_2)=[h_{1},h_{2}]_{\cH}$ for every $h_{1},h_{2}\in \cH$.
	\end{itemize}
\end{defn}

There are several classical examples for the choice of $\cH$. For instance, if $\cH:=L^{2}(0,T)$ then $W(t):=\cW 1_{[0,t]}$ defines a (real-valued) Brownian motion on $[0,T]$. Furthermore, if $\cH:=L^{2}(D)$ where $D\subseteq\mathbb{R}^d$ for $d\geq 1$ one obtains \emph{white noise} on $D$. 

\begin{defn}\label{zyl}
An $L^2(0,T;\cH)$-isonormal process is called {\bf $\cH$-cylindrical Brownian motion} on $[0,T]$.	
\end{defn}

An $\cH$-cylindrical Brownian motion will be denoted by $W_{\cH}$. From the definition one immediately has that:
\begin{itemize}
	\item $(W_{\cH}(t) h)_{t\in[0,T]}$ is a real-valued Brownian motion for every $h\in\cH$;
	\item $\mathbb{E} (W_{\cH}(t)h\cdot W_{\cH}(s)g)=\min\{s,t\}[g,h]_{\cH}$ for every $g,h\in\cH$ and $s,t\in[0,T]$.
\end{itemize}

Setting $\cH:=L^2(D)$ in the Definition~\ref{zyl} we remark that a $L^2(D)$-cylindrical Brownian motion provides the mathematical model for {\bf space-time white noise} on $[0,T]\times D$. This explains why $\cH$-cylindrical Brownian motions appear naturally in the context of stochastic partial differential equations.~As in the finite-dimensional case, one intuitively expects to represent $W_{\cH}$ as a series given by
\begin{align}\label{series:zyl}
W_{\cH}(t) = \sum\limits_{n=1}^{\infty} w_n(t)e_{n},
\end{align} 
where $(w_{n}(\cdot))_{n\geq 1}$ are one-dimensional independent standard Brownian motions and $(e_{n})_{n\geq 1}$ stands for an orthonormal basis in $\cH$. However, such a series does not always describe a genuine $\cH$-valued process, since it may not always converge in $L^2(\Omega,\mathcal{F},\mathbb{P};\cH)$. Therefore, one is interested in necessary conditions that ensure the convergence of the series~\eqref{series:zyl} in $L^{2}(\Omega,\mathcal{F},\mathbb{P};\cH)$.
One can show that there always exists a larger space $\overline{\cH}$ that contains $\cH$ such that~\eqref{series:zyl} converges in $L^{2}(\Omega,\mathcal{F},\mathbb{P};\overline{\cH})$, according to Proposition 2.5.2~in~\cite{DaPratoZabczyk} and Proposition 2.3.4~in~\cite{PrevotRoeckner}. Namely, if the embedding operator $\cH\hookrightarrow\overline{\cH}$ is Hilbert-Schmidt, i.e.~the map $\cH\ni x\mapsto x\in \overline{\cH}$ is Hilbert-Schmidt, then~\eqref{series:zyl} converges in $L^{2}(\Omega,\mathcal{F},\mathbb{P};\overline{\cH})$.
Conditions that ensure the convergence of~\eqref{series:zyl} in $\cH$ can be expressed in terms of trace-class operators on $\cH$. In this context we recall that a non-negative self-adjoint operator $Q\in\mathcal{L}(\cH)$ is called {\bf trace-class} if
\begin{align}\label{trace}
\mbox{Tr~} Q = \sum\limits_{n=1}^{\infty} [Q e_{n}, e_{n}]_{\cH} =\sum\limits_{n=1}^{\infty} \lambda_n <\infty,
\end{align}
where the sequence of eigenvectors $(e_{n})_{n\geq 1}$ of $Q$ forms an orthonormal basis of $\cH$ and $(\lambda_{n})_{n\geq 1}$ are the corresponding eigenvalues.
This leads to the following definition / characterization of an $\cH$-valued Wiener process. For more details, see Proposition~4.1~in~\cite{DaPratoZabczyk}.

\begin{defn}
	An $\cH$-valued stochastic process is called {\bf $Q$-Wiener process} (notation $(W^{Q}(t))_{t\in[0,T]}$) if the following two conditions hold:
	\begin{itemize}
		\item $(W^{Q}(t))_{t\in[0,T]}$ is a Gaussian process on $\cH$ with mean zero and covariance operator $t Q$ for $t\geq 0$;
		\item For any $t\geq 0$ the process $(W^{Q}(t))_{t\in[0,T]}$ has the following representation
		\begin{align}\label{series:q}
		W^{Q}(t) =\sum\limits_{n=1}^{\infty} \sqrt{\lambda}_{n} w_{n}(t) e_{n},
		\end{align}
		where $w_{n}(t):= \frac{1}{\sqrt{\lambda_{n}}}(W^{Q}(t), e_{n})$ for $n\in\mathbb{N}$ are independent real-valued standard Brownian motions on $(\Omega,\mathcal{F},\mathbb{P})$ and the series~\eqref{series:q} converges in $L^{2}(\Omega, \mathcal{F},\mathbb{P};\cH)$. As previously introduced, $(e_{n})_{n\geq 1}$ is the sequence of eigenvectors of $Q$ with corresponding eigenvalues $(\lambda_{n})_{n\geq 1}$.
	\end{itemize}
\end{defn}

The requirement that the trace of $Q$ must be finite can immediately be derived computing the variance of $W^Q$ as given by~\eqref{series:q} as follows:
\begin{align}
\mathbb{E} \|W^{Q}(t)\|^2_{\cH} = \mathbb{E} \Bigg(\sum\limits_{n=1}^{\infty} \lambda_{n} w^2_{n}(t) [e_{n},e_{n}]_{\cH}\Bigg)=\sum\limits_{n=1}^{\infty} \lambda_{n} \mathbb{E} w^{2}_{n}(t) = t\mbox{ Tr}~Q .
\end{align}
Of course, the covariance of $W^{Q}$ is given by
$$\mathbb{E}(W^{Q}(t)\cdot W^{Q}(s))=\min\{s,t\}\mbox{ Tr}~Q, ~ s,t\in[0,T]. $$
The construction of the stochastic integral of Hilbert-space valued processes with respect to a $Q$-Wiener process (more general cylindrical Wiener processes) can be seen as the natural generalization of the finite-dimensional It\^{o}-integral to Hilbert spaces. Further details can be looked up in~\cite[Section~4.2]{DaPratoZabczyk},~\cite{PrevotRoeckner,LiuRoeckner}. As in the finite-dimensional one has an analogue square-integrability condition together with an It\^o-isometry. 

\begin{thm}
	Let $\cU$ denote a further separable Hilbert space and let $\phi:[0,T]\times\Omega\to\mathcal{L}_2(\cU,\cH)$ be $\mathcal{B}([0,T]\otimes\mathcal{F};\mathcal{B}(\mathcal{L}_2(\cU,\cH))$-measurable such that 
	\begin{align*}
\mathbb{P} \Bigg(\int\limits_{0}^{T} \|\phi(s)\|^{2}_{\mathcal{L}_{2}(\cU,\cH)}~\txtd s<\infty\Bigg) =1.
	\end{align*}
	Then $\phi$ is stochastically integrable with respect to an $\cU$-cylindrical Brownian motion $W_{\cU}$ and the following {\bf It\^{o}-isometry} holds true:
	\begin{align}
	\mathbb{E} \Bigg\| \int\limits_{0}^{t} \phi(s)~\txtd W_{\cU}(s) \Bigg\|^2_{\cH} =\mathbb{E}\int\limits_{0}^{t} \|\phi(s)\|^2_{\mathcal{L}_{2}(\cU, \cH)}~\txtd s,~~\mbox{for  } t\in[0,T]. 
	\end{align}
\end{thm}

After introducing the preliminary notions on infinite-dimensional stochastic processes and integrals, we can now consider the SPDE
\begin{align}\label{spde1}
\begin{cases}
\txtd u(t) =  [A u(t) + f (t,u(t))]~\txtd t + g(t,u(t))~\txtd W_{\cU}(t),~~\mbox{} t\in[0,T]\\
u(0)=u_{0}\in \cH,
\end{cases}
\end{align}
where the assumptions on the coefficients are precisely specified in each section. As a typical example, the reader may already think of $u=(u(t))_{t\in[0,T]}=(u(t,x))_{t\in[0,T],x\in D}$ (where the space-dependency is usually dropped whenever there is no risk of confusion), the Laplacian $A=\Delta$ and sufficiently regular maps $f,g$. The aim of the remaining part of this section is to provide an overview and discuss various solution concepts for~\eqref{spde1}.

\subsection{Mild, Weak and Strong Solutions}\label{two:1}

Throughout this section the stochastic basis $(\Omega,\mathcal{F}, (\mathcal{F}_{t})_{t\in[0,T]}, \mathbb{P})$ is fixed and we make the following assumptions on the coefficients of~\eqref{spde1}. The initial condition $u_{0}$ is $\mathcal{F}_{0}$-measurable, the linear operator $A$ generates a $C_{0}$-semigroup $(S(t))_{t\geq 0}$ on $\mathcal{H}$, the nonlinear terms $f:[0,T]\times\Omega\times\mathcal{H}\to\mathcal{H}$, $g:[0,T]\times\Omega\times\cH \to \mathcal{L}_{2}(\cU,\cH)$ are $(\mathcal{B}([0,t])\otimes\mathcal{F}_{t}\otimes\mathcal{B}(\cH); \mathcal{B}(\cH) )$-measurable, respectively $(\mathcal{B}([0,t])\otimes\mathcal{F}_{t}\otimes\mathcal{B}(\cH); \mathcal{B}(\mathcal{L}_{2}(\cU,\cH)) )$-measurable for every $t\in[0,T]$ and $(W_{\cU}(t))_{t\in[0,T]}$ is an $\cU$-cylindrical Brownian motion. For notational simplicity, the $\omega$-dependence of $f$ and $g$ has been dropped. \\
Similar to the deterministic PDE case, we recall the following solution concepts for~\eqref{spde1}.

\begin{defn}\label{strong}
An $\cH$-valued adapted process $(u(t))_{t\in[0,T]}$ having $\mathbb{P}$-a.s. Bochner integrable trajectories is called a {\bf strong solution} for~\eqref{spde1} if for all $t\in[0,T]$ 
\begin{align*}
\int\limits_{0}^{t} u(s)~\txtd s \in \textnormal{Dom}(A),~~\mathbb{P}-a.s.,
\end{align*}
and for all $t\in[0,T]$:
\begin{align*}
u(t) = u_{0} + \int\limits_{0}^{t} [A u(s) + f(s,u(s))]~\txtd s  + \int\limits_{0}^{t} g(s,u(s))~\txtd W_{\cU}(s), ~~\mathbb{P}-a.s. 
\end{align*}
\end{defn}

\begin{defn}\label{weak:pde}
An $\cH$-valued adapted process $(u(t))_{t\in[0,T]}$ having $\mathbb{P}$-a.s. Bochner integrable trajectories is called a {\bf weak solution} for~\eqref{spde1} if for every test function $\zeta\in \textnormal{Dom}(A^{*})$ and $t\in[0,T]$ we have:
\begin{align*}
\left<u(t),\zeta\right> &=\left<u(t),u_{0}\right> + \int\limits_{0}^{t} [\left<u(s), A^{*}\zeta\right> + \left<f(s,u(s)),\zeta\right>] ~\txtd s\\
& + \int\limits_{0}^{t} \left<g(s,u(s)),\zeta\right>~\txtd W_{\cU} (s), ~~\mathbb{P}-a.s.
\end{align*}
\end{defn}

Clearly, a strong solution is also a weak one.

\begin{defn}\label{def:mild} An $\cH$-valued adapted process $(u(t))_{t\in[0,T]}$ is called {\bf mild solution} for~\eqref{spde1} if 
\begin{align}
\label{int:mild}
\mathbb{P} \Bigg( \int\limits_{0}^{T} \|u(s)\|^2_{\cH}~\txtd s<\infty\Bigg) =1
\end{align}
	and for $t\in[0,T]$ the {\bf variation of constants} (or {\bf Duhamel's}) {\bf formula} holds true:
	\begin{align}\label{mild:sol}
	u(t) &= S(t) u_{0} + \int\limits_{0}^{t} S(t-s)f(s,u(s))~\txtd s\\
	& + \int\limits_{0}^{t} S(t-s)g(s,u(s))~\txtd W_{\cU}(s),~~ \mathbb{P}-a.s.
	\end{align}
\end{defn}

The first integral in~\eqref{mild:sol} will be referred to as {\bf deterministic convolution} and the second one will be called {\bf stochastic convolution}. This is well-defined due to the assumption~\eqref{int:mild}. If additionally the $C_{0}$-semigroup $(S(t))_{t\geq 0}$ is {\emph analytic} (\cite[Chapter~3]{Pazy}), then one can derive optimal regularity results for the stochastic convolution~\cite[Sections~5.4 and~6.4]{DaPratoZabczyk},~\cite{NeervenVeraarWeis, NVW}.
Using a classical fixed-point argument, one can prove the following existence result of mild solutions for~\eqref{spde1}, see Theorem 7.2~in~\cite{DaPratoZabczyk}.

\begin{thm}\label{mild:ito}
	Let $f$ and $g$ additionally satisfy the following Lipschitz and growth boundedness assumptions:
\begin{itemize}
	\item [1)] there exists a constant $L>0$ such that for all $u,v\in\cH$, $t\in[0,T]$ and almost all $\omega\in\Omega$ we have:
	\begin{align}\label{lipschitz:mild}
	\|f(t,\omega,u)- f(t,\omega,v)\|_{\cH} + \|g(t,\omega,u) -g(t,\omega,v)\|_{\mathcal{L}_{2}(\cU, \cH)} \leq L \|u-v\|_{\cH};
	\end{align}
	\item [2)] there exists a constant $l>0$ such that for all $u\in\cH$, $t\in[0,T]$ and almost all $\omega\in\Omega$ we have:
	\begin{align}
		\|f(t,\omega,u)\|^2_{\cH} + \|g(t,\omega,u) \|^2_{\mathcal{L}_{2}(\cU,\cH)} \leq l^2 (1+\|u\|^2_{\cH}).
 	\end{align}
\end{itemize}
	 There exists a unique (up to equivalence) mild solution of~\eqref{spde1}.
\end{thm}

Naturally one can obtain similar assertions under weaker assumptions on the coefficients (e.g.~local Lipschitz continuity, dissipativity, local monotonicity etc.) using for instance cut-off and localization techniques, see~\cite{DaPratoZabczyk,PrevotRoeckner, LiuRoeckner}.

\begin{remark}
	\begin{itemize}
		\item [1)] 	Under the assumptions on the coefficients stated above, weak and mild solutions for~\eqref{spde1} are equivalent, see~Theorem~5.4~\cite{DaPratoZabczyk} and~\cite[Theorem~12]{Ondrejat}. A meaningful example where strong, weak and mild solutions are equivalent is given by the linear stochastic heat equation perturbed by a trace-class Brownian motion,~\cite[Theorem~5.14]{DaPratoZabczyk} and~\cite[Theorem~ 8.10]{isem}.
		\item [2)] The solution concepts discussed in this subsection are {\bf strong in the probabilistic sense}, meaning that the probability space together with the processes defined on them have been fixed.
	\end{itemize}
\end{remark}

We point out that this existence result of mild solutions for SPDEs carries over with suitable modifications to the Banach space-valued case, see~\cite{NeervenVeraarWeis, isem} and the references specified therein. First of all, one considers
Banach spaces which have certain geometric properties (which are satisfied by all relevant reflexive spaces appearing in the PDE theory such as Sobolev, Besov, Bessel potential spaces) and replaces the Hilbert-Schmidt operators by $\gamma$-radonifying ones (Definition~5.8~\cite{isem}).
Furthermore, one has to impose a stronger Lipschitz assumption in~\eqref{lipschitz:mild}, a so-called "randomized Lipschitz condition" (\cite[Section~5]{NeervenVeraarWeis}). This can be interpreted as a Gaussian version of Lipschitz continuity. This is mainly necessary due to the fact that if a Banach space $\cV$ has the property that $g(u)$ is stochastically integrable (with respect to a Brownian motion) for every $\cV$-valued function $u$ and Lipschitz continuous function $g:\cV\to\cV$, then $\cV$ is isomorphic to a Hilbert space. Therefore, it is highly non-trivial to find an appropriate Banach space in order to set up a fixed-point argument and conclude existence results for Banach space-valued evolution equations in the It\^o-setting. In Section~\ref{rough:paths} we present a pathwise approach, where such technical difficulties are not encountered.

\subsection{Martingale Solutions}\label{two:two}

The notion of martingale solution for an SPDE is the same as weak solution for a finite-dimensional SDE. In this context, since one also has weak in the PDE sense, recall Definition~\ref{weak:pde} one refers to weak solutions in the probabilistic case as \emph{martingale solutions}.

\begin{defn} If for given data $\cH$, $Q$, $u_{0}\in\cH$ and coefficients $A, f, g$ there exists a probability space
	$(\Omega,\mathcal{F},\mathbb{P})$ together with a filtration $(\mathcal{F}_{t})_{t\geq 0}$ and a $Q$-Wiener process such that $u$ is the mild solution of the SPDE
	\begin{align}\label{spde:m}
	\begin{cases}
	\txtd u =  [A u + f (u)]~\txtd t + g(u)~\txtd W^{Q}(t),~~\mbox{} t\in[0,T]\\
	u(0)=u_{0}\in \cH,
	\end{cases}
	\end{align}
	then $(\Omega,\mathcal{F},\mathbb{P}, (\mathcal{F}_t)_{t\geq 0}, W^{Q}, u)$ is called a {\bf martingale solution} for~\eqref{spde:m}.
\end{defn}

This means that for an initial configuration one seeks a probability space together with a $Q$-Wiener process which is defined on it, such that $u$ is a mild solution for~\eqref{spde:m}. This means that the probability space and $W^Q$ are part of the solution. The general strategy for existence of martingale solutions is based on compactness arguments,~\cite[Section 8.3]{DaPratoZabczyk} or on the Girsanov theorem, see~\cite[Section~10.3]{DaPratoZabczyk}.
We shortly describe the first situation which additionally assumes that the semigroup $(S(t))_{t\geq 0}$ is~\emph{compact} in order to apply Arzela Ascoli's theorem. For example, this holds true if $A$ is the realization of a uniformly elliptic operator on a bounded domain and it will also be explored in Chapter~\ref{dynamics}.
\begin{thm}
	Let $(S(t))_{t\geq 0}$ be compact, $f:\cH\to \cH$, $g:H\to\mathcal{L}_{2}(\cH)$ be continuous mappings and $W^{Q}$ a $Q$-Wiener process on $\cH$. If additionally $f$ and $g$ satisfy the linear growth condition
	\begin{align}\label{growth:m}
	\|f(u)\|_{\cH} + \|g(u)\|_{\mathcal{L}_{2}(\cH)} \leq c (1+ \|u\|_{\cH}), ~~x\in\cH,
	\end{align}
	then the SPDE~\eqref{spde:m} has a martingale solution.
\end{thm}
The proof of this theorem relies on the following steps.
\begin{itemize}
	\item [1)] One constructs solutions $u^n$ for regular coefficients $f^n$ and $g^n$ on a probability space $(\Omega,\mathcal{F},\mathbb{P})$ with a filtration $(\mathcal{F}_{t})_{t\geq 0}$ and a $Q$-Wiener process $W^{Q}$. Here $f^n:=P_{n}(f((P_{n}u)))$ and $g^n:=P_{n}(g((P_n u)))$ 
	for $u\in \cH$ where $P_{n}$ denotes the projection on the finite-dimensional subspace spanned by $\{e_{1},\ldots e_{n}\}$. Again $(e_{n})_{n\geq 1}$ stands for an orthonormal basis in $\cH$. The main idea is to approximate $f^n$ and $g^n$ by Lipschitz continuous mappings that satisfy~\eqref{growth:m}~(with possibly another constant $c$ independent on $n$).
	\item [2)] One shows that the sequence of laws $\{\mbox{Law}(u^{n})\}_{n\geq 1}$ converges weakly in $C([0,T],\cH)$ to a measure $\mu$.
	\item [3)] One constructs the solution $u$ of~\eqref{spde:m} with the law $\mu$ on a new probability space $(\tilde{\Omega},\mathcal{\tilde{F}},\mathbb{\tilde{P}})$ for a new filtration $(\mathcal{\tilde{F}}_{t})_{t\geq 0}$ with respect to a new $Q$-Wiener process $\tilde{W}^{Q}$. The main tools required to perform this step are Skorohod's embedding and martingale representation theorems.
\end{itemize}

For further details on this topic, see~\cite[Chapter~8]{DaPratoZabczyk} and~\cite{BrMo14}.

\subsection{Variational Solutions}
\label{variational}

In contrast to the previous sections, we  deal now with {\em nonlinear operators} $A$ satisfying suitable monotonicity, coercivity and growth conditions and discuss the variational approach of~\cite{LiuRoeckner, BarbuRoeckner1, RoecknerWang, BarbuPratoRoeckner, PrevotRoeckner}. Further solution concepts for nonlinear operators will be presented in Section~\ref{quasilinear}.
Similar to~\eqref{spde1} we consider
\begin{align}\label{spde2}
\begin{cases}
\txtd u(t) =  A (t,u(t)) ~\txtd t + g(t,u(t))~\txtd W_{\cU}(t),~~\mbox{} t\in[0,T]\\
u(0)=u_{0}\in \cH.
\end{cases}
\end{align}
Let $\cV$ stand for a reflexive Banach space such that the embeddings $\cV\hookrightarrow\cH\hookrightarrow \cV^{*}$ are continuous and dense. Then $(\cV,\cH, \cV^{*})$ is a {\bf Gelfand triple}. 
The coefficients
$$A: [0,T]\times \cV \times \Omega\to \cV^{*} \mbox { and } g:[0,T]\times \cV \times\Omega\to \mathcal{L}_{2}(\cU,\cH) $$
satisfy suitable measurability conditions and the following assumptions:
\begin{itemize}
	\item [1)] (Hemicontinuity) For all $u,v,x\in \cV$, $\omega\in\Omega$ and $t\in[0,T]$ the map
	\begin{align*}
	\mathbb{R}\ni \lambda \mapsto \hspace*{2 mm}_{\cV^{*}}\langle A(t,u+\lambda v,\omega),x\rangle_{\cV}
	\end{align*}
	is continuous.
	\item [2)] (Local monotonicity) There exists $c\in\mathbb{R}$ such that for all $u,v\in\cV$
	\begin{align*}
	2 \hspace*{1 mm}_{\cV^{*}}\langle A(\cdot,u) -A(\cdot,v), u-v\rangle_{\cV} &+ \|g(\cdot,u) - g(\cdot,v)\|^{2}_{\mathcal{L}_{2}(\cU,\cH)} \\
	&\leq c \|u-v\|^{2}_{\cH} \mbox{ on  } [0,T]\times\Omega.
	\end{align*}
	\item [3)] (Coercivity) There exists $\alpha\in(1,\infty)$ and constants $c_{1}\in\mathbb{R}$ and $c_{2}\in(0,\infty)$ and an $(\mathcal{F}_{t})$-adapted process $\overline{f}\in L^{1}([0,T]\times\Omega, \txtd t \otimes \mathbb{P})$ such that for all $v\in\cV$ and $t\in[0,T]$
	\begin{align}\label{coercivity}
	2 \hspace*{1 mm}_{\cV^{*}}\langle A(t,v), v\rangle_{\cV} +\|B(t,v)\|^{2}_{\mathcal{L}_{2}(\cU,\cH)} \leq c_{1} \|v\|^{2}_{\cH} - c_{2}\|v\|^{p}_{\cV} + \overline{f}(t), ~\mbox{on  } \Omega.
	\end{align}
	\item [4)] (Boundedness) There exist $c_{3}\in[0,\infty)$ and an $(\mathcal{F}_{t})$-adapted process $h\in L^{\frac{p}{p-1}}([0,T]\times\Omega, \txtd t \otimes\mathbb{P})$ such that for all $v\in\cV$ and $t\in[0,T]$
	\begin{align*}
\|A(t,v)\|_{\cV^{*}} \leq h(t) + c_{3}\|v\|^{p-1}_{\cV}, ~~\mbox{on  } \Omega,
	\end{align*}
	where $p$ is the constant from assumption 3).
\end{itemize}

\begin{ex}
	Let $D\subset \mathbb{R}$ be an open bounded domain and consider the Gelfand triple \begin{align*}
	\cV:=W^{1,2}_{0}(D)\hookrightarrow L^{2}(D) \hookrightarrow W^{-1,2}(D).
	\end{align*}
	We assume for simplicity that $g:\cV\to\mathcal{L}_{2}(L^{2}(D))$ is Lipschitz continuous and point out the following operators:
	\begin{itemize}
		\item $Au:=\Delta u + f(u)\nabla u$, where $f:\mathbb{R}\to\mathbb{R}$ is a bounded Lipschitz function,
		\item and $Au:= \Delta u - |u|^{m-2}u +\eta u$, for $\eta<0$ and $1\leq m<2$,
	\end{itemize}
which both satisfy the properties introduced above.
\myendex 

For further applications, see~\cite[Section~3]{LiuRoeckner},~\cite[Section~3]{GessLiuRoeckner},~\cite{BarbuPratoRoeckner,BarbuRoeckner1, PrevotRoeckner} and the references specified therein. Additional examples are also provided in Section~\ref{stability}.

\begin{defn}
A continuous $\cH$-valued process $u$ is called a {\bf variational solution} for~\eqref{spde2} if for its $\txtd t\otimes\mathbb{P}$-equivalent class $\overline{u}$ we have the following:
	\begin{itemize}
		\item $\overline{u}\in L^{p}([0,T]\times\Omega,\txtd t\otimes\mathbb{P}, \cV)\cap L^{2}([0,T]\times\Omega,\txtd t\otimes \mathbb{P},\cH)$,
		where $p$ is the constant from assumption 3).
		\item  
		$u(t)=u_{0} + \int\limits_{0}^{t} A(s,\overline{u}(s))~\txtd s + \int\limits_{0}^{t} g(s,\overline{u}(s))~\txtd W_{\cU}(s), ~~\mathbb{P}-a.s., t\in[0,T]$.\\
	\end{itemize}
\end{defn} 
The main existence and uniqueness result in this setting is given in~\cite[Theorem~1.1]{LiuRoeckner}.
\begin{thm}
	Under the above assumptions, for any given initial condition $u_{0}\in L^{2}(\Omega,\mathcal{F}_{0},\mathbb{P})$ there exists a unique variational solution of~\eqref{spde2}.
\end{thm}
The proof of this statement relies on a Galerkin approximation combined with suitable a-priori estimates of the solution. The variational approach is also applicable if one further perturbs~\eqref{spde2} by a general additive fractional noise as in~\cite{RoecknerWang}. Results regarding L\'evy-type noise can be looked up in~\cite{GessLiuRoeckner}. Results regarding the equivalence of mild and variational solutions can be found in~\cite[Section~7.5]{Veraar}.
 
\subsection{The Quasilinear Case}\label{quasilinear}

Using similar notations as in Section~\ref{variational} we are now interested in {\em semigroup methods} and {\em mild solutions} for SPDEs constituted by
\begin{align}\label{spde3}
\begin{cases}
\txtd u (t) =  A (u(t))u(t) ~\txtd t + g(u(t))~\txtd W_{\cU}(t),~~\mbox{} t\in[0,T]\\
u(0)=u_{0}\in \cH.
\end{cases}
\end{align}
To this aim we let for simplicity $g: \cH\to\mathcal{L}_2(\cU,\cH)$ be Lipschitz continuous. Furthermore we denote by $\cP_{T}$ the set of all $\cH$-valued continuous $(\mathcal{F}_{t})$-adapted processes on $[0,T]$ and impose the following restrictions on $A$:
\begin{itemize}
	\item [1)] For every $u\in\cP_{T}$, $A(u)$ is a {\bf sectorial operator} of angle $0<\phi<\frac{\pi}{2}$, namely  $$\sigma(A(u))\subset \Sigma_{\phi}:=\{ \lambda\in\mathbb{C}: 
	|\mbox{arg}~ \lambda|<\phi \}, \mbox{ for } u\in \mathcal{P}_{T},  $$
	where $\sigma(A(u))$ denotes the spectrum of $A(u)$.
	\item [2)] For $u\in \mathcal{P}_{T}$ the resolvent operator 
	$(\lambda \mbox{Id} -A(u))^{-1}$ satisfies the {\bf Hille-Yosida estimate}, i.e., there exists 
	$\widetilde{M}\geq 1$ such that
	$$||(\lambda \Id -A(u))^{-1}||_{\mathcal{L}(\cH)} \leq \frac{\widetilde{M}} {|\lambda|+1},\mbox{ 
		for } \lambda\notin\Sigma_{\phi} \mbox{ and } u\in \mathcal{P}_{T}.$$
	\item [3)] Let $0<\nu\leq 1$ be fixed and $L\geq 1$. Then
	\begin{equation}\label{nu}
	||A^{\nu}(u) (A(u)^{-1}- A(v)^{-1}) ||_{\mathcal{L}(\cH)}\leq 
	L || u-v||_{\cH}, \mbox{ for } u,v\in \mathcal{P}_{T}.
	\end{equation}
\end{itemize}

The first two assumptions guarantee that for $u\in\cP_{T}$, $A(u)$ generates a {\bf parabolic evolution system} $U^{u}$, which is a family of linear operators 
depending on two-time parameters and having similar properties with analytic 
$C_{0}$-semigroups, consult~\cite[Section~II.2]{Amann4} and~\cite[Chapter~3]{Pazy}. Naturally, since $u$ is a stochastic process, $U^{u}$ is obviously also $\omega$-dependent. One has for every $\omega\in\Omega$ that
\begin{align*}
& U^{u}(t,t,\omega)=\mbox{Id},~~\mbox{ for all } t\in[0,T];\\
&  U^{u}(t,s,\omega)U^{u}(s,r,\omega)=U^{u}(t,r,\omega),~~\mbox{  for all } r\leq s\leq t.
\end{align*}
Assumption 3) represents a classical Lipschitz continuity, which can be weakened according to~\cite{KuehnNeamtu}. One can easily incorporate a locally Lipschitz drift term. Keeping this in mind, one would expect to obtain by fixed-point methods a mild solution for~\eqref{spde3}. This should be given by the variation of constants formula
\begin{align}\label{mild:false}
u(t)= U^{u}(t,0,\omega) u_{0} + \int\limits_{0}^{t} U^{u}(t,s,\omega)g(u(s))~\txtd W_{\cU}(s).
\end{align}
However, this cannot hold true, since the terms appearing in the stochastic convolution in~\eqref{mild:false} do not satisfy the necessary properties required in order to define the It\^o-integral. More precisely, one can show that~\cite[Proposition~2.4]{Veraar} the mapping $\omega\mapsto U^{u}(t,s,\omega)$ is only $\mathcal{F}_{t}$-measurable but {\em not $\mathcal{F}_s$-measurable}. A possible ansatz to overcome this fact is to use the Skorohod integral for non-adapted integrands and Malliavin calculus~\cite{Nualart}.
Another possibility is to exploit the integration by parts formula in order to define~\eqref{mild:false} in a meaningful way. This concept of solution is referred to {\em pathwise mild} and was developed in~\cite{PronkVeraar} for {\em nonautonomous linear random} operators $A$. 

\begin{defn}
	A continuous $\cH$-valued adapted process $(u(t))_{t\in[0,T]}$ is called a {\bf pathwise mild solution} for~\eqref{spde3} if it satisfies
	\begin{align}\label{pm}
	u(t)&= U^{u}(t,0,\omega)u_{0} + U^{u}(t,0,\omega) \int\limits_{0}^{t} g(u(s))~\txtd W_{\cU}(s) \nonumber\\
	&- \int\limits_{0}^{t} A(u(s))U^{u}(t,s,\omega)\int\limits_{s}^{t} g(u(r))~\txtd  W_{\cU}(r)~\txtd s,~~\mathbb{P}-a.s.
	\end{align}
\end{defn}
The last two terms in the previous formula form  a {\em generalized stochastic convolution}.  The term~\eqref{pm} is well-defined due to the $\mathbb{P}$-a.s.~H\"older continuity of the integral $\int\limits_{0}^{\cdot}g(u(s))~\txtd W_{\cU}(s) $ which is necessary in order to compensate the fact that $\|A(u(s))U^{u}(t,s,\omega)\|_{\mathcal{L}(\cH)}\leq C(\omega)(t-s)^{-1}$ for $0\leq s<t\leq T$. In the nonautonomous case, one can also use the {\em forward integral} of Russo-Vallois~\cite{RussoVallois} to define the convolution~\eqref{mild:false} and to construct a solution to the corresponding SPDE which coincides with the pathwise mild solution as argued in~\cite[Section~4.5]{PronkVeraar}. 
\begin{thm}
	Under the previous assumptions, the quasilinear problem~\eqref{spde3} has a unique local-in-time pathwise mild solution.
\end{thm}
The statement can be proved using fixed-point arguments as in~\cite{KuehnNeamtu}. This can be achieved in two steps. First of all, fixing a process $v\in \cP_{T}$ and setting $A_{v}(t):=A(v(t,\omega))$, we obtain a Cauchy problem with time-dependent random drift
\begin{align}\label{v}
\txtd u = A_{v}(t) u(t)~\txtd t + g(u(t))~\txtd W_{\cU}(t),~~~ t\in[0,T].
\end{align}
Using \cite[Theorem 5.3]{PronkVeraar} we infer that~\eqref{v} has a pathwise 
mild solution given by
\begin{align*}
u(t) & = U^{v}(t,0)u_{0} + U^{v}(t,0)\int\limits_{0}^{t}g(u(s))~\txtd W_{\cU}(s) \\
& - \int\limits_{0}^{t}A_{v}(s)U^{v}(t,s)A_{v}(s) \int\limits_{s}^{t}
g(u(r)) ~\txtd W_{\cU}(r)~\txtd s,  ~  \mathbb{P}-a.s.,
\end{align*}
where $U^{v}(t,s)$ is the random parabolic evolution operator generated 
by $A_{v}$. For the sake of brevity we dropped the $\omega$-dependence in the previous formula.
The second step is to prove that the mapping
$$\Phi(v):=u , \mbox{ for } v\in \mathcal{P}_{T}$$ maps $\cP_{T}$ into itself and is a contraction if the time-horizon $T$ is chosen sufficiently small.
As applications we mention the stochastic Shigesada-Kawasaki-Teramoto model from mathematical biology, which was introduced to analyze population segregation by induced cross-diffusion, see~\cite[Section~4]{KuehnNeamtu}.  Written in divergence form, this is given by
\begin{align*}
\txtd u = \mbox{div}(\widetilde{\cA}(u)\nabla u)~\txtd t + f (u)~\txtd t + g(u)~\txtd W_{\cU}(t),
\end{align*}
where for $u:=(u_{1},u_{2})^{T}$, the matrix $\widetilde{\cA}$ is constituted by
$$\widetilde{\cA}(u):=\begin{pmatrix}
k_{1}+  2 c u_{1}+ a u_{2} & a u_{1} \\
b u_{2} & k_{2}+2 d u_{2} +b u_{1},
\end{pmatrix}
$$
and the nonlinear term is the same as in the classical Lotka-Volterra model, i.e.~
$$
f(u):=\begin{pmatrix}
\delta_{11} u_{1} - \gamma_{11} u_{1}^{2} -\gamma_{12}u_{1}u_{2}\\
\delta_{21} u_{2} - \delta_{21} u_{1}u_{2} -\gamma_{22} u_{2}^{2}
\end{pmatrix}.$$
The constants above are assumed to be positive and chosen such that the matrix $\widetilde{\cA}(u)$ is positive definite.\\
We conclude this subsection pointing out that there are numerous other solution concepts for quasilinear equations such as martingale solutions, strong (in the probabilistic sense) and weak (in the PDE sense) solutions~\cite{HofmanovaZhang}, where the coefficients are approximated with locally monotone ones, as discussed in Section~\ref{variational}, kinetic solutions~\cite{DHV,FehrmanGess}, entropy solutions~\cite{GK}, solution concepts for quasilinear SPDEs via rough paths theory~\cite{OttoWeber}, paracontrolled calculus~\cite{FurlanGubinelli, BDH} and regularity structures~\cite{GerencserHairer}.

\subsection{Random Field Approach}
\label{random:field}

In this section we briefly present the main ideas of an alternative method to~\cite{DaPratoZabczyk} (Section~\ref{two:1}) to construct solutions for SDPEs, see~\cite{Dalang,DalangLluis} for further details. As described at the beginning of Section~\ref{chapter2}, the framework of~\cite{DaPratoZabczyk} relies on defining stochastic integrals with respect to Hilbert space-valued processes, whereas the random field approach is based on Walsh integration theory~\cite{Walsh} with respect to martingale measures. This approach will also be explored in Section~\ref{dynamics} in order to compute moments of the solutions and Lyapunov-exponents. As a reminder, the concept {\bf random field} refers to a family of random variables indexed by several parameters, usually space and time. We now fix $d\geq 1$ and consider the heat equation
\begin{align}\label{heat:randomfields}
\frac{\partial{u(t,x)}}{\partial t}= \Delta u(t,x) + u(t,x)\dot{W}(t,x),~~t>0 \mbox{ and } x\in\mathbb{R}^{d},
\end{align}
where $\dot{W}=\partial_{t} W$ and the noise is white in time and (possibly) correlated in space,
i.e.~ $$\mathbb{E}[\dot{W}(t,x)\dot{W}(s,y)]=\delta_{0}(t-s)f(x-y),~~s,t>0\mbox{ and  }x,y\in\mathbb{R}^d.$$
Here $\delta_{0}$ is the Dirac distribution. We further assume that the initial condition $u_{0}$ is constant and $f$ is a positive definite function having a {\bf spectral measure} $\mu$, i.e., $\mathcal{F}\mu=f$, where $\mathcal{F}$ denotes the Fourier-transform. There are several examples for $f$ such as:
\begin{itemize}
	\item [1)] $f=\delta_{0}$. In this case one obtains {\bf space-time white-noise} (on $\mathbb{R}_{+}\times\mathbb{R}^{d}$), recall Definition~\ref{zyl}.The spectral measure is $\mu(\txtd \xi)= (2\pi)^{-d/2}\txtd\xi$.
	\item [2)] ({\bf Riesz kernels}) $f(x)=|x|^{-\eta}$ for $0<\eta<d$. Here the spectral measure is $\mu(\txtd\xi)=c_{d,\eta}|\xi|^{-(d-\eta)}\txtd\xi$, where the constant $c_{d,\eta}$ depends only on the dimension $d$ and on $\eta$.
	\item [3)] ({\bf Fractional kernels}) For $d=1$ and $H\in(0,1)$ let $f(x)=c_{H}|x|^{2H-2}$ with $\mu(\txtd\xi)=c_{H}|\xi|^{1-2H}\txtd\xi$. We notice that if $H\leq 1/2$, $f$ is not a measure anymore.
\end{itemize}
Let $\overline{\cU}$ stand for a separable Hilbert space which will be described below. For $G:[0,T]\times\Omega\to\overline{\cU}$ (satisfying suitable measurability assumptions), we introduce 
the notation
\begin{align}\label{walsh}
\int\limits_{0}^{T} G(s)~\txtd W(s):=\int\limits_{0}^T \int\limits_{\mathbb{R}^{d}} G(s,y) W(\txtd s, \txtd y)
\end{align}
and emphasize that the stochastic integral~\eqref{walsh} is defined in the sense of Walsh~\cite{Walsh}. More precisely, letting $\overline{\cU}$ be the Hilbert space obtained by the completion of $C^{\infty}_0(\mathbb{R}^d)$ under 
\begin{align*}
\|\phi\|^{2}_{\overline{\cU}}= (2\pi)^{d/2} \int\limits_{\mathbb{R}^d} \mu(\txtd\xi) ~|\mathcal{F}(\phi)(\xi)|^2,
\end{align*}
one has for an adapted process $G\in L^{2}([0,T]\times\Omega,\overline{\cU})$ the {\em isometry}
\begin{align}\label{isometry:rf}
\mathbb{E} \Bigg( \int\limits_{0}^{T} G(s)~\txtd W(s)\Bigg)^2 =\mathbb{E} \|G\|^{2}_{L^2(0,T,\overline{\cU})} = \mathbb{E} \int\limits_{0}^{T} \txtd t \int\limits_{\mathbb{R}^{d}} \mu(\txtd\xi) ~|\mathcal{F}(G(t))(\xi)|^{2}.
\end{align}

\begin{defn}
	We call a stochastic process $\{u(t,x) : t>0, x\in\mathbb{R}^d\}$ satisfying
	\begin{align}\label{msol:rf}
	u(t,x)= \Gamma_{t} u_{0} (x) + \int\limits_{0}^{t}\int\limits_{\mathbb{R}^d} \Gamma_{t-s}(y-x) u(s,y) W(\txtd s,\txtd y),~~\mathbb{P}-a.s.,
	\end{align}
	for all $t>0$ and $x\in\mathbb{R}^d$, a {\bf mild random field solution} of~\eqref{heat:randomfields}.
\end{defn}
Here $\Gamma$ is the fundamental solution of the heat equation / heat kernel, i.e.~
\begin{align*}
\Gamma_{t}(x) =(4\pi t)^{-d/2} \exp\Big(\frac{-|x|^2}{4t}\Big)
\end{align*}
and the convolution appearing in~\eqref{msol:rf} is constructed via Walsh's theory.\\
In this case, one can easily compute the Fourier-transform of the heat kernel, i.e.~
\begin{align*} 
\mathcal{F}(\Gamma_{t})(\xi) =\exp(-4\pi^2t|\xi|^2),~~\xi\in\mathbb{R}^d
\end{align*}
and observe that
\begin{align*}
\int\limits_{0}^{T} \exp(-4\pi^2t|\xi|^2) ~\txtd t =\frac{1}{4\pi^2|\xi|^2} (1-\exp(-4\pi^2T|\xi|^2)).
\end{align*}
Recalling~\eqref{isometry:rf}, in order that
\begin{align*}
\int\limits_{0}^{T} \txtd t \int\limits_{\mathbb{R}^{d}} \mu(\txtd\xi) ~|\mathcal{F}(\Gamma_t)(\xi)|^{2}<\infty,
\end{align*}
one has to impose that the spectral measure satisfies
\begin{align}\label{spectral:measure}
	\int\limits_{\mathbb{R}^d} \frac{\mu(\txtd \xi)}{1+|\xi|^2} <\infty.
\end{align}
Condition~\eqref{spectral:measure} is referred to as {\bf Dalang's condition}. For space-time white-noise, recall $\mu(\txtd\xi)=(2\pi)^{-d/2}\txtd\xi$,  it is immediately clear that~\eqref{spectral:measure} holds true only in $d=1$. Consequently, in higher space-dimensions, the problem~\eqref{heat:randomfields} becomes ill-posed.~This aspect was investigated in~\cite{Hairer} and will also be addressed in Section~\ref{reg:structures}.\\
Existence of mild random field solutions for~\eqref{heat:randomfields} was obtained in~\cite{Dalang}.
In this framework, one has an analogue statement to Theorem~\ref{mild:ito}. Namely, if one incorporates measurable drift and diffusion coefficients $f$ and $g$ in~\eqref{heat:randomfields} (using the same notations as in Theorem~\ref{mild:ito}) satisfying Lipschitz and growth boundedness conditions, one obtains a mild random field solution 
$\{u(t,x)\mbox{ : } t>0, x\in \mathbb{R}^{d}\}$ such that for every $p\geq 1$
	\begin{align*}
	\sup\limits_{t\geq 0}\sup\limits_{x\in\mathbb{R}^d}\mathbb{E}|u(t,x)|^{p}<\infty.
	\end{align*}
The proof uses Picard iteration, see~\cite{Dalang},~\cite[Theorem~4.3]{DalangLluis}.


\subsection{Mild Solutions for Rough SPDEs}\label{rough:paths}

In this section we would like to deal with a random forcing, which unlike Brownian motion is no longer a semi-martingale and its trajectories are not independent. A meaningful example in this sense is given by the {\bf fractional Brownian motion} (fBm), originally introduced by B.~Mandelbrot and J.~van Ness in~\cite{MandelbrotNess}. This is also a Gaussian process, but its covariance function additionally depends on a Hurst index / parameter $H\in(0,1)$. If $H=1/2$ one obtains Brownian motion, if $H\neq 1/2$, this process exhibits different behavior from Brownian motion. The key feature which is required to comprehend the next steps is the H\"older continuity of its trajectories for an exponent $\alpha<H$. The entire solution theory relies only on this regularity result and is developed within the framework of the {\em rough path approach}~\cite{Lyons, Gubinelli, GubinelliTinde, FritzHairer}.
Similar to Skorohod integrals~\cite{Nualart}, rough path techniques also go beyond the semi-martingale case, where It\^o calculus and Walsh theory no longer apply. The major difference is that these provide a completely pathwise construction of stochastic integrals and implicitly of the solution of SPDEs. However, many solution concepts and techniques that are available in the context of deterministic PDEs have not been completely developed for rough SPDEs (e.g.~weak solutions, variational inequalities, etc). Therefore, this has been a very active research field in the last few years with numerous recent developments, see for instance~\cite{DejaM, HHofmanova, Gubinelli, GubinelliLejayTindel, DeyaGubinelliTinde, HairerG} and the references specified therein.\\
In the finite-dimensional case, the theory is well-understood~\cite{FritzHairer, FritzVictoir, Lyons} and rough path techniques have been successfully applied to investigate dynamical aspects of solutions, see~\cite{BailleulRiedelScheutzow, KN, RiedelDelay}.\\
In this section, we provide the main intuition of the construction of a mild solution for a SPDE similar to~\eqref{spde1}, but driven by multiplicative fractional noise. For further details regarding the following approach, see~\cite{HesseNeamtu1}. 
We fix the {\bf Hurst parameter} $H\in(1/3,1/2]$ and consider the non-linear SPDE~\eqref{spde:m}, driven by a trace-class fBm $B^{H}$, more precisely 
\begin{align}\label{spde:fbm}
\begin{cases}
\txtd u = [A u + f(u)]~\txtd t + g(u)~\txtd B^{H}(t), ~~ t\in[0,T]\\
u(0)=u_{0}\in\mathcal{H}.
\end{cases}
\end{align}
Here, $B^{H}$ is a $\cU$-valued fBm and is defined as in~\eqref{series:q} taking a sequence of independent standard real-valued fBm-s $(w_n)_{n\in\mathbb{N}}$ having the same Hurst index $H\in(1/3,1/2]$.
Furthermore, we assume that $A$ generates an \emph{analytic} $C_{0}$-semigroup and for $\beta\in\mathbb{R}$ we denote by $\mathcal{H}_{\beta}:=\textnormal{Dom}((-A)^{\beta})$ the domains of its fractional powers. These can be identified with Sobolev spaces for certain ranges of $\beta$. For a detailed description of these spaces, see~\cite[Chapter~3]{Pazy}. 
The diffusion coefficient $g:\mathcal{H}\to\mathcal{L}(\cU,\cH)$ is supposed to be three times continuously differentiable with bounded derivatives and also Lipschitz continuous as a mapping $g:\cH\to\mathcal{L}(\cU, \cH_{\beta})$. Typical examples of such operators are given by integral operators with smooth kernels, see~\cite{HesseNeamtu1,GarridoLuSchmalfuss, DeyaGubinelliTinde}. Naturally, for finite-dimensional noise one can consider polynomials with smooth coefficients, see~\cite{DeyaGubinelliTinde, HairerG}. For simplicity we set $f=0$, since this term does not require any additional arguments.
Similar to~\eqref{def:mild} we call a mild solution for~\eqref{spde:fbm} a process $u$ that satisfies the variation of constants formula
\begin{align}\label{mild:fbm}
u(t)= S(t)u_{0}  + \int\limits_{0}^{t} S(t-s)g(u(s))~\txtd B^{H}(s).
\end{align}
However, additional measurability and adaptedness conditions are not required. We now provide a purely pathwise definition of the stochastic convolution
\begin{align}\label{integral}
\int\limits_{0}^{t} S(t-s) g(u(s))~\txtd B^{H}(s).
\end{align}
The strategy to define~\eqref{integral} relies on an approximation procedure. We firstly consider a smooth path $B^{H}$ and a (H\"older) continuous trajectory $u$. The general argument eventually follows considering smooth approximations of $B^{H}$. Thereafter, the passage to the limit entails a suitable construction/interpretation of all the expressions above according to~\cite[Section~5]{HesseNeamtu1}. \\Throughout this section we use the standard rough path notation, i.e.~the value at time $t$ of a function $y$ is given by $y_{t}$ instead of $y(t)$. For $y:[0,T]\to \cH$, $y_{st}:=y_{t}-y_{s}$ denotes an {\bf increment}. We use the same notation to indicate the evaluation at two time points $s$ and $t$ of an element $y:\Delta_{T}\to\overline{\cH}$, where $\overline{\cH}$ stands for an arbitrary separable Hilbert space. For the semigroup $S$ we keep the notation $S(\cdot)$.\\
Our aim is to define~\eqref{integral} using Riemann-Stieltjes sums and a Taylor expansion for $g$. In the following $\mathcal{P}$ stands for a partition of the interval $[0,t]$. By a {\em formal computation}, this ansatz reads as

\begin{align}
&\hspace*{-30 mm}\int\limits_{0}^{t}{S(t-r)g(u_r)}~\txtd B^{H}_r 
= \sum\limits_{[v_{1},v_{2}]\in \mathcal{P}} S(t-v_{2}) \int\limits_{v_{1}}^{v_{2}}{S(v_{2}-r)g(u_r)}~\txtd B^{H}_r \nonumber\\
\approx \sum\limits_{[v_{1},v_{2}]\in \mathcal{P}} S(t-v_{2}) \Big[ &\int\limits_{v_{1}}^{v_{2}}{S(v_{2}-r)g(u_{v_{1}})} ~\txtd B^{H}_{r}\nonumber \\
&+ \int\limits_{v_{1}}^{v_{2}}{S(v_{2}-r)\txtD g(u_{v_{1}})(u_{r}-u_{v_{1}})}~\txtd B^{H}_r\Big] \nonumber\\
=: &\sum\limits_{[v_{1},v_{2}]\in \mathcal{P}} S(t-v_{2}) 
\Big[\omega_{v_{1}v_{2}}^S(g(u_{v_{1}})) +z_{v_{1}v_{2}}(\txtD g(u_{v_{1}}))\Big]. \label{heuristic:integral}
\end{align}
Here we introduced the notation
\begin{align}\label{omegas_h}
\omega^{S}_{v_{1}v_{2}}(G(y_{v_{1}})):=\int\limits_{v_{1}}^{v_{2}} S(v_{2}-r) g (u_{v_{1}}) ~\txtd B^{H}_{r},
\end{align}
respectively

\begin{align}\label{z_h}
z_{v_{1}v_{2}}(\txtD g(y_{v_{1}})):= \int\limits_{u}^{v} S(v_{2}-r)  \txtD g(u_{v_{1}}) (u_{r} - u_{v_{1}}) ~\txtd B^{H}_{r}.
\end{align}
By a classical integration by parts formula, see Theorem~3.5~in~\cite{Pazy} one can argue that the term $\omega^{S}$ can be defined for a rough input $B^{H}$. This is no longer the case for $z$. By a classical rough path ansatz, where one {\em plugs in the definition of the solution itself in the construction of the rough integral}, one obtains the following approximation for $z$:
\begin{align*}
&z_{st}(E)  \approx \sum\limits_{[v_{1},v_{2}]\in \mathcal{P}} S(t-v_{2}) \Big[ b_{v_{1}v_{2}}(E,g(u_{v_{1}})) + a_{v_{1}v_{2}}(E,u_{v_{1}})\Big] - \omega^S_{st}(E u_{s}),
\end{align*}
where $E$ is a placeholder that stands for $\txtD g$ and
\begin{align*}
b_{v_{1}v_{2}}(E,g(u_{v_{1}})) :=\int\limits_{v_{1}}^{v_{2}} S(v_{2}-r) E \int\limits_{v_{1}}^{r} S(r-q) g (u_{v_{1}}) ~\txtd B^{H}_{q}~\txtd B^{H}_{r},
\end{align*}
respectively
\begin{align*}
a_{v_{1}v_{2}}(E,u_{v_{1}}):=\int\limits_{v_{1}}^{v_{2}} S(v_{2}-r) E S (r-v_{1}) 
u_{v_{1}} ~\txtd B^{H}_{r}.
\end{align*}
This indicates that we have to define $a$, $b$ 
and $\omega^{S}$ in order to fully characterize $z$.
As already emphasized, the previous formal computation has been conducted under the assumption that $B^{H}$ is a smooth path. In this case $a, b$ and $\omega^{S}$ are well-defined. However, the main challenge is to construct these processes for rough inputs $B^{H}$ and to show that~\eqref{heuristic:integral} is indeed the right way to define~\eqref{mild:fbm}. These arguments contain suitable approximation techniques (\cite[Section~5]{HesseNeamtu1}) combined with the additional necessary assumption $g:\cH\to\mathcal{L}(\cU, \cH_{\beta})$, since it is not possible to define the {\em iterated integral} $b$ for arbitrary operators~\cite[Remark~4.3]{DeyaGubinelliTinde}.
Another essential tool is represented by the Sewing Lemma, which ensures the existence of the rough integral as a limit of Riemann-Stieltjes sums, see~\cite[Theorem~4.1]{HesseNeamtu1} and the references specified therein. Putting all these together, one concludes that the {\em pathwise} solution of~\eqref{spde:fbm} is given by a pair $(u,z)$, where
\begin{align}
 u_{t}&= S(t)u_{0} + \sum\limits_{[v_{1},v_{2}]\in \mathcal{P}} S(t-v_{2}) 
\Big[\omega_{v_{1}v_{2}}^S(g(u_{v_{1}})) +z_{v_{1}v_{2}}(\txtD g(u_{v_{1}}))\Big]\label{u}\\
 z_{st}(Dg) &= \sum\limits_{[v_{1},v_{2}]\in \mathcal{P}} S(t-v_{2}) \Big[ b_{v_{1}v_{2}}(\txtD g(u_{v_{1}}),g(u_{v_{1}})) + a_{v_{1}v_{2}}(\txtD g(u_{v_{1}}),u_{v_{1}})\Big]\nonumber\\
 & \hspace*{13 mm}- \omega^S_{st}(\txtD g(u_s))\label{z}.
\end{align}
Here $u$ is referred to the {\em path} component of the solution and $z$ is the {\em area term}.
\begin{thm}
	The rough SPDE~\eqref{spde:fbm} has a unique local-in-time solution $(u,z)$, where the two components $u$ and $z$ are given by~\eqref{u} and~\eqref{z}.
\end{thm}

The existence of such a solution is obtained by a fixed-point argument in a suitable function space, which is chosen in order to compensate the time-singularity in zero of $(S(t))_{t\geq 0}$. This is required since the terms appearing in the Riemann-Stieltjes approximations must exhibit a certain H\"older regularity to make the Sewing Lemma applicable. However, using estimates for analytic semigroups, one has for a fixed $\gamma\in (0,1]$ that
\begin{align}
\|S(t) u_{0} - S(s) u_{0}\|_{\cH} &=\|(S(t-s) -\mbox{Id}) S(s) u_{0}\|_{\cH}\nonumber\\
& \leq C (t-s)^{\gamma} \|u_{0}\|_{\cH_{\gamma}}\nonumber\\
& \leq C (t-s)^{\gamma} s^{-\gamma}\|u_{0}\|_{\cH}. \label{h:s}
\end{align}
Keeping~\eqref{h:s} in mind, one is motivated to introduce for $\gamma_1, \gamma_2\in(0,1]$ the function space
\begin{align*}
C^{\gamma_{1},\gamma_{2}}([0,T],\cH):=\Bigg\{y\in C([0,T],\cH) ~:~\sup\limits_{0<s<t\leq T} s^{\gamma_{1}} \frac{\|y_{t} -y_{s}\|_{\cH}}{(t-s)^{\gamma_{2}}}<\infty\Bigg\}.
\end{align*}
Regarding the area component, we also have
\begin{align*}
C^{\gamma_1,\gamma_2}(\Delta_{T},\overline{\cH}):=\Bigg\{y\in C(\Delta_{T},\overline{\cH}) ~:~\sup\limits_{0<s<t\leq T} s^{\gamma_{1}} \frac{\|y_{st} \|_{\overline{\cH}}}{(t-s)^{\gamma_{2}}}<\infty\Bigg\},
\end{align*}
for a certain Hilbert space $\overline{\cH}$.
Using~\eqref{h:s} one has that $S(\cdot)u_{0}\notin C^{\gamma}([0,T],\cH)$ but $S(\cdot)u_{0}\in C^{\gamma,\gamma}([0,T];\cH)$.
Regarding this aspect together with further necessary regularity conditions required for the Sewing Lemma, the appropriate space for the fixed-point argument turns out to be
\begin{align*}
\Bigg\{(u,z)~:~ u\in  C^{\beta,\beta}([0,T],\cH) \mbox{ and }\\ z\in C^{\alpha}(\Delta_{T},\mathcal{L}(\mathcal{L}(\cH\otimes\cU,\cH) ,\cH)) \times C^{\alpha+\beta,\beta}(\Delta_{T},\mathcal{L}(\mathcal{L}(\cH\otimes\cU,\cH) ,\cH)) \Bigg\} .
\end{align*}
Here $\alpha<H$ stands for the H\"older regularity of $B^{H}$ and $\beta$ is chosen such that $\alpha+2\beta>1$ and gives the time-regularity of the solution of~\eqref{spde:fbm}.\\
Due to the Taylor expansion employed in~\eqref{heuristic:integral} one obtains certain quadratic estimates of the solution, therefore the global-in-time existence is not straightforward. This follows e.g.~under an additional boundedness assumption on $g$, using regularizing properties of analytic semigroups, as argued in~\cite{HesseNeamtu2, DeyaGubinelliTinde}.\\
Before pointing out some concluding remarks, we emphasize that for finite-dimensional noise, the solution theory for~\eqref{spde:fbm} simplifies. In this case, one can show that the solution $(u,z)$ is given by
\begin{align}\label{mild:fbm:finitedim}
(u,z)= \Bigg(S(\cdot)u_{0}+\int\limits_{0}^{\cdot} S(\cdot-r)g(u_{r})~\txtd B^{H}_r, g(u)\Bigg),
\end{align}
similar to the ODE case~\cite[Chapter~8]{FritzHairer}. This can be achieved by means of a fixed-point argument in a function space that incorporates additional {\em space regularity} in order to compensate the missing time-regularity in~\eqref{h:s}. More precisely, one can prove that the following approximation 
\begin{align*}
\int\limits_{0}^{t} S(t-r)g(u_{r})~\txtd B^{H}_{r} &=\lim\limits_{|\mathcal{P}|\to 0} \Bigg(\sum\limits_{[v_{1},v_{2}]\in\mathcal{P}} S(t-v_{1}) g (u_{v_{1}}) B^{H}_{v_{1}v_{2}} \\
&\hspace*{-4 mm} + \sum\limits_{[v_{1},v_{2}]\in\mathcal{P}}S(t-v_{1}) \txtD g(u_{v_{1}})g(u_{v_{1}}) \int\limits_{v_{1}}^{v_{2}} (B^{H}_{r} - B^{H}_{v_{1}})~\txtd B^{H}_{r} \Bigg)
\end{align*}
can be used to define the stochastic convolution.
 For further details regarding this topic, see~\cite{HairerG}.

\begin{remark}
	\begin{itemize}
		\item [1)] The techniques are applicable also in the Banach space-valued setting. For finite-dimensional noise, the solution of~\eqref{spde:fbm} is given by~\eqref{mild:fbm:finitedim}.
		For infinite-dimensional noise, one  needs to ensure that the series~\eqref{series:q} defines a Banach space-valued fractional Brownian motion, see~\cite[Section~3]{BzNS} and~\cite[Proposition 8.8]{isem}.
		\item [2)] We emphasize that rough paths techniques are not restricted to the case of fBm, but apply to a wider class of {\em Gaussian processes}, their covariance functions satisfy certain criteria, as specified in~\cite[Chapter~10]{FritzHairer} and~\cite[Chapter~15]{FritzVictoir}.
		\item [3)] If the trajectories of the fBm are more regular, i.e.~$H\in(1/2,1)$, then~\eqref{mild:fbm} can be defined using the Young integral, see~\cite{GubinelliLejayTindel}, as
		\begin{align*}
		\int\limits_{0}^{t} S(t-r) g(u_{r})~\txtd B^{H}_{r} =\lim\limits_{|\mathcal{P}|\to 0} \sum\limits_{[v_{1},v_{2}]\in\mathcal{P}} S(t-v_{2}) \omega^{S}_{v_{1}v_{2}} (g(u_{v_{1}})).
		\end{align*}
		\item [4)] Another possibility to introduce mild solutions for such SPDEs is to decompose the integral in a series of one-dimensional integrals and define these using fractional calculus, see~\cite{GarridoLuSchmalfuss, MaslowskiNualart}. However, such an approach leads to restrictions on the coefficients and on the trace of the noise. Using this ansatz, one needs to impose as in~\cite{MaslowskiNualart} that $\mbox{Tr}~Q^{1/2}=\sum\limits_{n=1}^{\infty}\sqrt{\lambda_{n}}<\infty$, compare~\eqref{trace}. Another choice would be to take $g:\cH\to\mathcal{L}_{2}(\cU,\cH)$ as in~\cite{GarridoLuSchmalfuss}.
		\item [4)] If one considers~\eqref{spde:fbm} driven by additive fractional noise, then~\eqref{integral} can easily be defined as a suitable It\^o-type integral, see~\cite{DuncanMaslowski, MaslowskiNeerven}.
	\end{itemize}
\end{remark}

\subsection{Renormalized Solutions}\label{reg:structures}

Here we review briefly an approach for singular SPDEs using regularity structures~\cite{Hairer,Hairer2}. To illustrate the main issue, let us consider the stochastic bistable {\bf Nagumo equation} on the torus with space-time white noise, i.e., 
\be
\label{eq:ACwhite}
\partial_t u = \Delta u + f(u) +\xi,\qquad f(u):=u(1-u)(u-p),
\ee
where $u=u(t,x)$, $x\in\R^d/\Z^d=:\T^d$ for $d=2$ or $d=3$, $\xi=\xi(t,x)=\partial_t W^{\Id}(t,x)$ is space-time white noise, $\Delta$ is the Laplacian, and $p\in\R$ is a parameter; we remark that the case of a cubic nonlinearity also occurs in the {\bf Allen-Cahn/Chaffee-Infante/$\Phi^4$/Real-Ginzburg-Landau/Schl\"ogl}\\ models, which share similar features regarding existence theory of solutions. The fundamental issue of~\eqref{eq:ACwhite} is that the roughness of the space-time white noise forcing $\xi$ interplays with the nonlinearity to prevent the application of more classical solution concepts. To understand this, let us introduce the {\bf parabolic scaling} $\fs:=(2,1,1)$ for $d=2$ and $\fs:=(2,1,1,1)$ for $d=3$, which can be used to define a norm 
\benn
\|(t,x)\|_\fs:=|t|^{1/2}+\sum_{j=1}^d|x_j|,
\eenn
which is better adapted to the natural scaling of the operator $L:=\partial_t-\Delta$. Let $C^\alpha_\fs$ denote the space of H\"older continuous functions $\phi:\R\times \R^d \ra \R$ with the scaled norm, which is straightforward to define for $\alpha\geq 0$. For $\alpha<0$, fix $r=-\lfloor \alpha \rfloor$ and let $\mathcal{B}^r_{\fs,0}$ denote the space of all $C^r$-functions supported on $\{z\in\R^{d+1}:\|z\|_\fs \leq 1\}$. A Schwartz distribution $v\in \mathcal{S}(\R^{d+1})$ belongs to $C^\alpha_\fs$ if it belongs to the dual space of $C_0^r$ and for every compact set $\mathfrak{K}\subset \R^{d+1}$ there exists a constant $C>0$ such that 
\benn
\langle v,\cS_{\fs,z}^\delta \eta\rangle \leq C\delta^\alpha ,\quad \forall\delta\in(0,1],~\forall z\in\mathfrak{K},~\forall \eta \in \mathcal{B}^r_{\fs,0} \text{ with }\|\eta\|_{C^r}\leq 1,
\eenn 
where $(\cS_{\fs,z}^\delta \eta)(\bar{t},\bar{z}):=\delta^{-(d+2)}\eta(\delta^{-2}(t-\bar{t}),\delta^{-1}(x-\bar{x}))$. Plainly, $C^\alpha_\fs$ generalized $\alpha$-H\"older functions do not blow up worse near each point than a power divergence with exponent $\alpha$. One can prove that space-time white noise satisfies
\be
\label{eq:stwhite}
\xi\in C^\alpha_\fs \quad \text{for $\alpha=-\frac{d+2}{2}-\kappa$ and any $\kappa>0$.}
\ee
Suppose we would want to solve~\eqref{eq:ACwhite} by re-formulating it as a fixed-point of an iterated map $M$
\be
u^{(k+1)}=L^{-1}(f(u^{(k)})+\xi)=:M(u^{(k)}).
\ee
Suppose we start with $u^{(0)}=0$, then we get $M(0)=L^{-1}(0+\xi)=L^{-1}(\xi)=u^{(1)}$. Now $u^{(1)}$ has higher regularity in comparison to $\xi$ since applying $L^{-1}$ corresponds to convolution with the heat kernel. Classical Schauder theory tells us that $u^{(1)}\in C^{\alpha+2}_\fs$ so we gained two orders of regularity. For $d=2$ we get $\alpha+2=-2-\kappa+2=-\kappa<0$ while for $d=3$ we have $\alpha+2=-\frac52-\kappa+2=-\frac12-\kappa<0$. So in both cases, we get that $u^{(1)}$ is just a Schwartz distribution with negative H\"older regularity. However, if we now want to compute $u^{(2)}$, then we have to make sense of $f(u^{(1)})$, which is just not possible within classical distribution/generalized function theory as we can only multiply distributions, which have non-negative sums of their H\"older regularity exponents. A natural attempt to still use a fixed point approach is to first smooth the noise, e.g., via a mollifier
\benn
\rho_\varepsilon(t,x):=\varepsilon^{-(d+2)}\rho(t\varepsilon^{-2},x\varepsilon^{-1}),\qquad \xi_\varepsilon:=\rho_\varepsilon\ast \xi,
\eenn
where $\rho$ is a compactly supported function of integral $1$, and $\ast$ denotes the space-time convolution. Then it is classical that 
\be
\label{eq:ACwhite1}
\partial_t u_\varepsilon = \Delta u_\varepsilon + f(u_\varepsilon) +\xi_\varepsilon,
\ee
has smooth, even global-in-time, solutions for sufficiently regular initial. However, it turns out that directly taking a limit $\varepsilon\ra 0$ leads to divergence so the SPDE is singular. Classical considerations in phase transitions / bifurcation theory suggest that the divergence can be prevented if the equation is {\bf renormalized} to 
\be
\label{eq:ACwhite2}
\partial_t u_\varepsilon = \Delta u_\varepsilon + f(u_\varepsilon)+C(\varepsilon)u_\varepsilon +\xi_\varepsilon,
\ee
where $C(\varepsilon)=\cO(\ln\varepsilon)$ for $d=2$ and $C(\varepsilon)=\cO(\varepsilon^{-1})$ for $d=3$ are computable diverging functions as $\varepsilon \ra 0$. If the initial condition $u(0,\cdot)$ belongs to a sufficiently regular H\"older space, e.g., continuity is sufficient, then one can prove~\cite{Hairer} that there exists a sequence of local-in-time solutions $u_\varepsilon$ such that
\benn
u_\varepsilon\ra u_0 \quad \text{as $\varepsilon \ra 0$ in probability independently of $\rho$.} 
\eenn
The last statement is one example application of the theory of regularity structures, i.e., it establishes the existence of a {\bf renormalized solution} $u_0$. To prove the convergence, a very elaborate construction is necessary. The basic steps are (I) building a regularity structure and a model space adapted to the SPDE, (II) solving the fixed point problem in the abstract space of modeled distribution, (III) proving that the renormalized modeled distributions converge and one may reconstruct an actual distribution as a solution. Basically, this procedure has been worked out as an abstract result applicable to many singular SPDEs~\cite{BrunedChandraChevyrevHairer,BrunedHairerZambotti,ChandraHairer,Hairer}. We only illustrate part of the first step to discuss, to which classes of singular SPDEs, we may apply the theory of regularity structures.

The first step is to build a {\bf regularity structure} $(\cA,\cT,\cG)$, where $\cA\subset \R$ is an index set bounded below without accumulation points in $\R$, $\cT=\oplus_{\alpha\in\cA}\cT_\alpha$ is the {\bf model space} consisting of a graded sum of Banach spaces with $\cT_0$ isomorphic to $\R$, and $\cG$ is the structure group. We shall not discuss the structure group but try to illustrate, how $\cA$ and $\cT$ are constructed. The idea is to construct an abstract jet space of symbols capable of representing suitable regularity classes of functions. For $\cT_\alpha$ with $\alpha\in\N$, we simply take the space of homogeneous polynomials of degree $\alpha$ in $(d+1)$ variables $X_0,X_1,\ldots,X_d$ where $X_0$ represents the time variable, and $X_j$ is to be interpreted as a formal symbol. We count {\bf homogeneity} $|\cdot|_\fs$ of a polynomial in the parabolic scaling, i.e., 
\benn
|X^k|_\fs=|k|_s:=2k_0+\sum_{j=1}^d k_j,\quad X^k=X_0^{k_0}\cdots X_d^{k_d}.
\eenn      
We also set $|{\bf 1}|_\fs=0$ with the unit element ${\bf 1}$ spanning $\cT_0$. Furthermore, we let $\Xi$ be a symbol representing the noise with $|\Xi|_\fs=\alpha_0:=-(d+2)/2-\kappa$ for any fixed (small) $\kappa>0$. The convolution/integration against the heat kernel is represented by a map $\cI:\cT\ra \cT$ and $|\cI(\tau)|_\fs=|\tau|_\fs+2$ for any $\tau\in\cT$. Consider the abstract iterated map  
\benn
U^{(k+1)}=\cI(f(U^{(k)})+\Xi),
\eenn 
where we hope that iteration yields a fixed point. If we iterate the map, we get new symbols, e.g., with $U^{(0)}=0$ we get $U^{(1)}=\cI(\Xi)$ with $|\cI(\Xi)|_\fs=-\alpha_0+2$. Now one can iterate again, which yields second and third powers $\cI(\Xi)^2$ and $\cI(\Xi)^3$ as well as $\cI(\cI(\Xi))$ and so on including various combinations involving polynomials if they are represented in the initial condition. Now suppose we continue this procedure for all possible iterates and all possible symbols with non-negative homogeneity as initial condition. Then we check homogeneities and only keep symbols with negative homogeneity. For $d=2$ this yields
\benn
\Xi, ~\cI(\Xi)^3,~\cI(\Xi)^2,~\cI(\Xi), 
\eenn
with homogeneities $-2-\kappa$, $-3\kappa$, $-2\kappa$ and $-\kappa$ respectively. Carrying out the same calculation for $d=3$ is a good exercise yielding more symbols. For both cases, the procedure terminates only giving a finite number of negatively homogeneous symbols, so we can set $\cA$ as only containing all the negative homogeneities and use the negatively homogenous symbols to span $\cT_\alpha$ for $\alpha<0$. The condition that $\cA$ is bounded from below effectively means that the SPDE is {\bf locally subcritical} and the theory of regularity structures is designed for this class of singular SPDEs. Roughly speaking, for the class of locally subcritical singular SPDEs, the theory of regularity structures works and provides a local-in-time renormalized solution. The full procedure can be found in the works~\cite{BrunedChandraChevyrevHairer,BrunedHairerZambotti,ChandraHairer,Hairer}, while some particular 
examples can be found in~\cite{BerglundKuehn1,BerglundKuehn2,Hairer3,HairerShen, HairerG}. \\

Another breakthrough in the theory of singular SPDEs is constituted by the {\bf paracontrolled calculus} developed in~\cite{GIP}. This was successfully applied to the KPZ equation~\cite{GubinelliPerkowski}, the $\Phi^4_3$ model~\cite{CC} and singular quasilinear problems~\cite{FurlanGubinelli, BDH}. 
A comparison between the two methods can be looked up in~\cite{BH}.

\section{Dynamics: Concepts \& Results}\label{dynamics}

The aim of this section is to provide an overview and discuss recent developments regarding the long-time behavior of the solutions of SPDEs, again with a strong focus on stochastic reaction-diffusion equations.
Similar to the deterministic setting~\cite{KuehnBook1,Temam,SellYou}, one is interested to investigate the stability of steady states, predict if changes of stability (bifurcations) occur, analyze if trajectories become (exponentially) close to each other or look for sets where trajectories accumulate. There are several ways to describe such kind of phenomena in stochastic dynamics. We focus on the random dynamical system theory~\cite{Arnold}, which is used to define invariant sets such as manifolds or attractors for the flow generated by an SPDE. We also indicate several other possibilities to quantify the long-time behavior of SPDEs, such as invariant measures, large deviations or sample path approaches, depending also on the different noise terms driving these SPDEs, as discussed in Section~\ref{chapter2}.


\subsection{Random Dynamical Systems}
\label{rds}

As seen in Section~\ref{chapter2}, there are several different ways to perturb a PDE by a random input, such as trace-class (fractional) Brownian motion, space-time white noise, noise which is white in time and correlated in space, etc. Therefore, we introduce the next concept, which indicates a model of the noise driving the PDE. Recall that $\cH$ denotes an arbitrary separable Hilbert space and $(\Omega,\mathcal{F},\mathbb{P})$ stands for a probability space, which will further be specified later on.

\begin{defn}({\mbox Metric dynamical system (MDS)})
Let $\theta:\mathbb{R}\times\Omega\rightarrow\Omega$ be a family of $\mathbb{P}$-preserving transformations (meaning that $\theta_{t}\mathbb{P}=\mathbb{P}$ for all $t\in\mathbb{R}$) with the following properties:
	\begin{itemize}
		\item[(1)] the mapping $(t,\omega)\mapsto\theta_{t}\omega$ is $(\mathcal{B}(\mathbb{R})\otimes\mathcal{F},\mathcal{F})$-measurable;
		\item[(2)] $\theta_{0}=\mbox{Id}_{\Omega}$;
		\item[(3)] $\theta_{t+s}=\theta_{t}\circ\theta_{s}$ for all $t,s,\in\mathbb{R}$.
	\end{itemize}
	Then the quadruple $(\Omega,\mathcal{F},\mathbb{P},(\theta_{t})_{t\in\mathbb{R}})$ is called a {\bf metric dynamical system}.
\end{defn}

In order to simplify the notation we write $\theta_{t}\omega$ for $\theta(t,\omega)$. We always assume that $\mathbb{P}$ is {\bf ergodic} with respect to $(\theta_{t})_{t\in\mathbb{R}}$, which means that any invariant subset has zero or full measure. Next, we introduce now the MDS corresponding to a genuine $\cH$-valued process having stationary increments. As examples, one can consider a trace-class (fractional) Brownian motion as given in~\eqref{series:q}.

\begin{ex}\label{mds1}	
Let $C_{0}(\mathbb{R},\cH)$ denote the set of continuous $\cH$-valued functions, which are zero at zero and are equipped with the compact open topology. Furthermore, by taking $\mathbb{P}$ as the Wiener measure on $\mathcal{B}(C_{0}(\mathbb{R},\cH))$ having a trace-class covariance operator $Q$ on $\cH$ and applying Kolmogorov's theorem about the existence of a  continuous version yields the canonical probability space $\Omega:=(C_{0}(\mathbb{R},\cH),\mathcal{B}(C_{0}(\mathbb{R},\cH)),\mathbb{P})$. To obtain an ergodic metric dynamical system we introduce the {\bf Wiener shift}, which is defined by
	\begin{align}\label{shift}
	\theta_{t}\omega(\cdot{})=\omega(t+\cdot{})-\omega(t),\mbox{  for  } \omega\in C_{0}(\mathbb{R},\cH).
	\end{align}
Whenever we work with the RDS approach, i.e.~Sections~\ref{im} and~\ref{ra}, we only consider the subset of all $\omega\in C_{0}(\mathbb{R},\cH)$ with subexponential growth and work with this new metric dynamical system using the same notations $(\Omega,\mathcal{F},(\theta_{t})_{t\in\mathbb{R}}, \mathbb{P})$.
\myendex 

\begin{remark} Let $\overline{W}_{\cH}$ stand for a process with stationary increments having a trace-class covariance operator $Q$ on $\cH$. Then one can embed $\overline{W}_{H}$ into a canonical probability space, as given in Example~\ref{mds1} and identify
	 \begin{align*}
	\overline{W}_{\cH}(t,\omega) =\omega(t), \mbox{ for } \omega\in C_{0}(\mathbb{R},\cH).
	\end{align*}
	We will use this notation in Sections~\ref{ra} and~\ref{im}.
\end{remark}

In this framework, one often has to deal with random constants whose values have to be controlled along the orbits of $\theta$. Therefore, the following concept plays a key role. For more details, see Proposition~4.1.3~in~\cite{Arnold}.

\begin{defn}\label{tempered}
A positive real-valued random variable $Y$ on a MDS $(\Omega,\mathcal{F},\mathbb{P},(\theta_{t})_{t\in\mathbb{R}})$ is called {\bf tempered} if there exists a $(\theta_{t})_{t\in\mathbb{R}}$-invariant set of full measure such that 
$$
	\lim\limits_{t\to\pm\infty}\frac{\log^{+}Y(\theta_{t}\omega)}{|t|}=0.
$$
\end{defn}

\begin{remark}
Note that temperedness is equivalent to {\bf subexponential growth}. The only alternative to this situation in the ergodic case is that
		$$ \limsup\limits_{t\to\pm\infty}\frac{\log^{+}Y(\theta_{t}\omega)}{|t|}=\infty.$$
\end{remark}

In the deterministic theory of dynamical systems, one obtains under certain assumptions on the coefficients the flow property of the solution operator. Since our PDE is now perturbed by a non-autonomous, irregular forcing, it is natural to expect that the time-evolution of the noise must play a role in this flow property. We describe the precise mathematical formalism of this fact~\cite{Arnold}.

\begin{defn}({\mbox Random dynamical system (RDS)})
\label{def:rds}
A continuous {\bf random dynamical system} on $\cH$ over a metric dynamical system\\ $(\Omega,\mathcal{F},\mathbb{P},(\theta_{t})_{t\in\mathbb{R}})$ is a mapping $$\varphi:\mathbb{R^{+}}\times\Omega\times \cH\to \cH,\mbox{  } (t,\omega,x)\mapsto \varphi(t,\omega,x), $$ 
which is $(\mathcal{B}(\mathbb{R}^{+})\otimes\mathcal{F}\otimes\mathcal{B}(\cH),\mathcal{B}(\cH))$-measurable and satisfies:
	\begin{itemize}
		\item[1)] $\varphi(0,\omega,\cdot{})=\mbox{Id}_{\cH}$ for all $\omega\in\Omega$;
		\item[2)]$ \varphi(t+\tau,\omega,x)=\varphi(t,\theta_{\tau}\omega,\varphi(\tau,\omega,x)), \mbox{ } x\in \cH, \mbox{ } t,\tau\in\mathbb{R}^{+}\mbox{ and }\omega\in\Omega;$
		\item[3)] $\varphi(t,\omega,\cdot{}):\cH\to \cH$ is continuous for all $t\in\mathbb{R}^{+}$.
		\end{itemize}
\end{defn}

The property 2) is referred to as {\bf cocycle property}. If one drops the $\omega$-dependence, one obtains exactly the flow property from the deterministic case. Here we see that the shift $\theta_{t}\omega$ describes the evolution of the noise.
Keeping this mind, we observe that RDS can be interpreted as the generalization of non-autonomous deterministic dynamical systems.
Referring to~\cite{Arnold}, it is well-known that an It\^{o}-type stochastic ordinary differential equation generates a random dynamical system under natural assumptions on the coefficients.~This fact is based on the flow property, see~\cite{Kunita, Scheutzow}, which can be obtained by Kolmogorov's theorem about the existence of a (H\"older)-continuous random field with finite-dimensional parameter range, i.e.~the parameters of this random field are the time and the non-random initial data. However, the generation of a RDS from an SPDE (as considered in Sections~\ref{two:1}--\ref{random:field}) has been a long-standing open problem, since Kolmogorov's theorem breaks down for random fields parametrized by infinite-dimensional Hilbert spaces, see~\cite{Mohammed}. As a consequence it is not trivial how to obtain a RDS from an SPDE, since its solution is defined almost surely, which contradicts Definition~\ref{def:rds}, where all properties must hold for all $\omega\in\Omega$. Particularly, this means that there are exceptional sets which depend on certain parameters of the SPDE, and it is not clear how to define a RDS if more than countably many exceptional sets occur. This problem was fully solved only under very restrictive assumptions on the structure of the noise driving the SPDE. More precisely, for {\em additive} and {\em linear multiplicative} noise, one can perform certain {\em Doss-Sussmann-type transformations}~\cite{Doss, Sussmann} and reduce the corresponding SPDE to a PDE with random coefficients. Such a PDE can be solved pathwise for every realization of the noise. We provide here two examples of such transformations. Based on these we analyze attractors and invariant manifolds for SPDEs in Sections~\ref{im} and~\ref{ra}. 

\begin{ex}(\mbox{SPDE with additive noise})\label{transf1}
	Consider the SPDE on $\cH$
	\begin{align}\label{spde:transf1}
		\txtd u(t) = [A u(t) + f (u(t))]~\txtd t + \txtd W^Q(t),
	\end{align}
	where $W^{Q}$ is a trace-class Brownian motion in $\cH$. We further assume that the semigroup generated by $A$ is exponentially stable, i.e.~there exist constants $\widetilde{M}\geq 1$ and $\overline{\mu}>0$ such that $\|S(t)\|_{\cH}\leq \widetilde{M} \txte^{-\overline{\mu} t}\text{ for }t>0$. The unique {\em stationary} solution of the {\em linear} SPDE
	\begin{align}\label{lin:spde}
	\txtd u(t) = A u(t)~\txtd t +\txtd W^Q(t)
	\end{align}
	is given by the {\bf Ornstein-Uhlenbeck process}, which can be represented by the convolution
	\begin{align*}
\int\limits_{-\infty}^{t} S(t-s)~\txtd W^{Q}(s).
	\end{align*}
	Taking the canonical probability space corresponding to $W^{Q}$, as introduced in Example~\ref{mds1}, and identifying $W^{Q}(t,\omega)=\omega(t)$, for $\omega\in\Omega$ we introduce the process $(t,\omega)\mapsto Z(\theta_{t}\omega)$ as
	\begin{align}\label{ou}
	Z(\theta_{t}\omega) :=\int\limits_{-\infty}^{t} S(t-s)~\txtd \omega(s) =\int\limits_{-\infty}^{0} S(-s)~\txtd\theta_{t}\omega(s).
	\end{align}
	Subtracting~\eqref{ou} from~\eqref{spde:transf1} one obtains 
	\begin{align}\label{pde:transf1}
	\txtd u(t) = [A u(t) + f(u(t) + Z(\theta_{t}\omega )) ]~\txtd t,
	\end{align}
	which is a PDE with random nonautonomous coefficients.
\myendex 

\begin{ex}(\mbox{SPDE with linear multiplicative noise})
\label{transf2}
Let $W$ stand for a one-dimensional standard Brownian motion.~Consider the SPDE on $\cH$ with Stratonovich noise
	\begin{align}\label{spde:transf2}
		\txtd u(t) = [A u(t) + f (u(t))] ~\txtd t + u \circ \txtd W(t).
	\end{align}
Performing the transformation $u^{*}=u\txte^{-z(\omega)}$, where $z$ is the one-dimensional Ornstein-Uhlenbeck process and dropping the $*$-notation, one obtains
\begin{align}\label{pde:transf2}
\txtd u(t) = [A u(t) + z(\theta_{t} \omega ) u(t)]~\txtd t + \overline{f}(\theta_{t}\omega,u(t))~\txtd t,
\end{align}
where $\overline{f}(\omega,u):=\txte^{-z(\omega)}f(\txte^{z(\omega)}u)$.
\myendex 

In both of the previous cases, one can show that the solution operators of~\eqref{pde:transf1} and~\eqref{pde:transf2} generate random dynamical systems. Furthermore,
 the dynamical systems generated by the original SPDEs are equivalent with those generated by~\eqref{pde:transf1} and~\eqref{pde:transf2}. This means that it is enough to work with the transformed equations and transfer all the results such as fixed-points, manifolds or attractors to the initial SPDEs. However, for general nonlinear multiplicative noise, this technique is obviously no longer applicable.~As a consequence of this issue, dynamical aspects for SPDEs such as stability, Lyapunov exponents, multiplicative ergodic theorems, random attractors, random invariant manifolds have not been investigated in their full generality. We discuss situations in the following sections and point out how the techniques presented in Section~\ref{rough:paths} can be employed to investigate the dynamics of SPDEs driven by multiplicative noise. Further applications in this sense can be looked up in~\cite{BailleulRiedelScheutzow, DeyaStability, KN,Gao, GS,FehrmanGess} and the references specified therein.

\subsection{Stability}\label{stability}

In this section we discuss several stability concepts for SPDEs such as
\begin{itemize}
	\item stability in probability;
	\item almost sure exponential stability;
	\item moment exponential stability;
	\item metastability.
\end{itemize}

Common tools used to investigate such concepts rely on the existence of Lyapunov functionals~\cite{Khasminskii, Mao}, random dynamical systems methods, ergodic theory and Lyapunov exponents~\cite{Arnold1, Arnold, CaraballoDuanLuSchmalfuss,lian, LuNS, GNS, GS}, or large deviation theory~\cite{Freidlin3,FreidlinWentzell,BerglundGentzM}. There are numerous works dealing with stability statements for It\^o-type SPDEs. Arguments via Lyapunov functionals heavily use the Markov property of the solution and semi-martingale techniques. Therefore, it is a challenging open problem to get optimal stability results for SPDEs driven by fractional Brownian motion, such as~\eqref{spde:fbm}. Progress in this direction was made in~\cite{DeyaStability, GS, GNS}.\\
In this section we only refer to SPDEs driven by Brownian motion, where the precise assumptions are stated below. We write $u(t,u_{0})$ in order to refer to a solution $u$ of such an SPDE at time $t>0$ having $u_{0}\in\cH$ as initial condition. Again, the probability space $(\Omega,\mathcal{F},\mathbb{P})$ is fixed. Let $q\geq 1$.

\begin{defn}\label{stab}
	We call a \emph{global-in-time} solution $u$ of an SPDE:
	\begin{itemize}
		\item [1)] {\bf stochastically stable} or {\bf stable in probability} if for every pair $\varepsilon\in(0,1)$ and $r>0$, there exists a $\delta=\delta(\varepsilon,r)>0$ such that
		\begin{align*}
		\mathbb{P}(\|u(t,u_{0})\|_{\cH} < r\mbox{ for all } t\geq 0 ) \geq 1-\varepsilon,
		\end{align*}
		for $\|u_{0}\|_{\cH}<\delta$;
		\item [2)] {\bf almost sure (global) exponentially stable} if for all $u_{0}\in\cH$
		\begin{align}\label{ly:exp}
		\limsup\limits_{t\to\infty} \frac{1}{t} \log (\|u(t,u_{0})\|_{\cH}) < 0~ \mbox{ a.s.};
		\end{align}
		\item [3] {\bf $q$-th moment exponentially stable} if for all $u_{0}\in\cH$
		\begin{align*}
			\limsup\limits_{t\to\infty} \frac{1}{t} \log (\mathbb{E}\|u(t,u_{0})\|^q_{\cH}) < 0.
		\end{align*}
 	\end{itemize}
\end{defn}

In general moment and almost sure exponential stability do not imply each other, see~\cite{Mao, Khasminskii} for further details. For  linear SPDEs, the quantity on the left-hand side of~\eqref{ly:exp} is often called {\bf Lyapunov exponent}, see~\cite{Khasminskii}. For a related, yet different,  definition of Lyapunov exponents in random dynamical system theory, see~\cite[Section~3.2]{Arnold} and~\cite{Arnold1,FlandoliS,CaraballoDuanLuSchmalfuss, lian}.
\begin{remark}
	Note that one can also investigate {\bf local exponential stability}, meaning that there exists a {\em random neighbourhood} $\cN(\omega)$ of zero such that for $u_{0}\in \cN(\omega)$, the relation~\eqref{ly:exp} is satisfied, compare Definition~11~in~\cite{GNS}.
\end{remark}

We now illustrate the concepts introduced in Definition~\ref{stab} for general SPDEs satisfying the assumptions formulated in Section~\ref{variational}. We first provide abstract conditions which ensure the stability of solutions of SPDEs and present thereafter concise examples. Using the same notations as in Section~\ref{variational} we recall that $(\cV,\cH,\cV^{*})$ stands for a Gelfand triple and let $\overline{\beta}>0$ such that $\|x\|_{\cH}\leq \overline{\beta}\|x\|_{\cV}$.
We further impose the following {\em coercivity} condition, compare~\eqref{coercivity}. There exist constants $\overline{\alpha}>0$, $\overline{\mu}>0$, $\overline{\lambda}\in\mathbb{R}$ and a nonnegative continuous function $h$ such that
\begin{align}\label{c:stability}
2~ _{\cV^{*}}\langle A(t,v), v\rangle_{\cV} + \|g(t,v)\|^{2}_{\mathcal{L}_{2}(\cU,\cH)} \leq -\overline{\alpha} \|v\|^{p}_{\cV} + \overline{\lambda} \|v\|^2_{\cH} + h(t)\txte^{-\overline{\mu}t},~~v\in\cV, 
\end{align}
where $p>1$ and for arbitrary $\delta>0$ we have $\lim\limits_{t\to\infty}\frac{h(t)}{\txte^{\delta t}}=0$. Then the following results ensure the mean square and almost sure exponential stability for the solutions of SPDEs, compare~\cite[Theorems~2.2,2.3]{CaraballoLiu}.

\begin{thm}\label{mean:square:s}
Let assumptions 1), 2), 4) in Section~\ref{variational} and~\eqref{c:stability} hold true. Then, if $u$ is a global solution of~\eqref{spde2}, there exists constants $\varepsilon>0$, $C>0$ such that
\begin{align*}
\mathbb{E}\|u(t)\|_{\cH}^2\leq C \txte^{-\varepsilon t}, ~~ t\geq 0,
\end{align*}
if either one of the following hypotheses are satisfied:
\begin{itemize}
	\item [a)] $\overline{\lambda}<0$ (for every $p>1$);
	\item [b)] $\overline{\lambda}\cdot\overline{\beta}^2-\overline{\alpha}<0$ (for $p=2$).
\end{itemize}
\end{thm}

\begin{thm}\label{as:st}
Under the assumptions of Theorem~\ref{mean:square:s}, there exist positive constants $\widetilde{M}$, $\varepsilon$ and a subset $N_{0}\subset\Omega$ with $\mathbb{P}(N_{0})=0$ such that for each $\omega\not\in N_{0}$, there exists a positive random number $\widetilde{T}(\omega)$ such that
\begin{align*}
\|u(t)\|^2_{\cH} \leq \widetilde{M} e ^{-\varepsilon t}, ~~t\geq\widetilde{T}(\omega).
\end{align*}
\end{thm}
The proofs of the previous theorems rely on the It\^o-formula (\cite[Theorem 4.2.5]{PrevotRoeckner}), martingale arguments and the Gronwall Lemma. For a better comprehension we now provide examples where the previous abstract stability criteria are applicable.
\begin{ex}
Let $a\in\mathbb{R}$, $b:\mathbb{R}\to\mathbb{R}$ be a Lipschitz continuous function with $b(0)=0$ and let $W$ stand for a one-dimensional standard Brownian motion.	We consider a stochastic heat equation on $\cH:=L^{2}(0,\pi)$ with Dirichlet boundary conditions given by
	\begin{align}
	\begin{cases}
	\txtd u =[ \Delta u + a u ]~\txtd t + b(u)~\txtd W(t), ~~t>0, ~~x\in(0,\pi)\\
	u(t,0)=u(t,\pi)=0, t>0\\
	u(0,x)=u_{0}(x).
	\end{cases}
	\end{align}
	Here $\cV:=W^{1,2}_{0}(0,\pi)$, $A(t,u):=\Delta u + au$ and $g(t,u):=b(u)$. One can verify that for $u\in\cV$:
	\begin{align*}
	2~ _{\cV^{*}}\langle A(t,u),u\rangle_{\cV} +\|g(t,u)\|^{2}_{\mathcal{L}_{2}(\mathbb{R},\cH)} \leq -2\|u\|^{2}_{\cV} + (2 a + l^2)\|u\|_{\cH},
	\end{align*}
	where $l$ is the Lipschitz constant of $b$. Therefore~\eqref{c:stability} holds for $p=2$ and in order to apply Theorem~\ref{mean:square:s} with $\overline{\alpha}=2$ and $\overline{\lambda}=2a+l^2$, we have to choose $\overline{\beta}>0$ such that $(2a + l^2)\overline{\beta}^2-2<0$, which can be achieved e.g.~setting $\overline{\beta}:=\frac{\pi}{\sqrt{2}}$.	
\myendex 

For a nonlinear operator $A$, we provide the following example.

\begin{ex}
	Let $a>0$ and $b:\mathbb{R}\to\mathbb{R}$ be a Lipschitz continuous with Lipschitz constant $l>0$ and assume $b(0)=0$. Furthermore, $W$ stands for a one-dimensional standard Brownian motion. We let $2<p<\infty$ and consider the SPDE on $\cH:=L^{2}(0,1)$
	\begin{align*}
	\begin{cases}
	\txtd u(t,x) = \Big[\frac{\partial }{\partial x} \Big( \Big|\frac{\partial u(t,x)}{\partial x} \Big|^{p-2}\frac{\partial u(t,x)}{\partial x} \Big) -a u(t,x) \Big]~\txtd t + b (u(t,x))~\txtd W(t),\\
	u(t,0)=u(t,1)=0,~~ t>0\\
	u(0,x)=u_{0}(x),\hspace*{9 mm} x\in[0,1].
	\end{cases}
	\end{align*}
	Here, $\cV:=W^{1,p}_{0}(0,1)$, $g(t,u)=b(u)$ and the nonlinear monotone operator $A:\cV\to\cV^{*}$ is defined as
	\begin{align*}
	_{\cV^{*}}\langle Au,v \rangle_{\cV} =\int\limits_{0}^{1} \Big| \frac{\partial u(x)}{\partial x}\Big|^{p-2}\frac{\partial u}{\partial x} \frac{\partial v}{\partial x}~\txtd x  -  a\int\limits_{0}^{1} u(x) v(x)~\txtd x, ~~u,v\in\cV.
	\end{align*}
	One can check the~\eqref{c:stability} holds true for $h\equiv 0$, $p>2$  $\overline{\alpha}=2$ and $\overline{\lambda}=-\kappa$, where $\kappa>0$ is chosen such that $l^2<2a-\kappa$, see~\cite[Example~3.2]{CaraballoLiu}. In this case one obtains again due to Theorems~\ref{mean:square:s} and~\ref{as:st} the mean square and almost sure exponential stability of the solution.
\myendex 
A very simple example, where almost all sample paths of the solution {\em do not} tend exponentially to zero is constituted by the following SDE.
\begin{ex}
	Let $a$, $b>0$ be two constants and $W$ stand for a one-dimensional standard Brownian motion. We consider the SDE
	\begin{align}\label{sde}
	\begin{cases}
	\txtd u = -a u(t)~\txtd t + (1+t)^{-b}~\txtd W(t),~~t>0\\
	u(0)=0.
	\end{cases}
	\end{align}
Setting $A(t,x):=-a x$, $g(t,x):=(1+t)^{-b}$ and letting $\langle \cdot,\cdot \rangle$ denote the usual scalar product in $\mathbb{R}$, we easily have that
\begin{align*}
2~\langle A(t,x), x \rangle + \|g(t,x)\|^{2} =-2a x^2+(1+t)^{-2b},
\end{align*}
where the last term obviously does not exponentially tend to zero.  Therefore~\eqref{c:stability} is not satisfied. Here, using the law of iterated logarithm, one can immediately verify that for the solution of~\eqref{sde}
\begin{align*}
u(t)= \txte^{-a t} \int\limits_{0}^{t} \txte^{a s}(1+s)^{-b}~\txtd W(s),
\end{align*}
it holds that
\begin{align*}
\limsup\limits_{t\to\infty}\frac{1}{t} \log |u(t)| = 0, ~~\mbox{a.s.}
\end{align*}
Taking another diffusion coefficient, e.g.~$g(t,x):=\txte^{-b t}$, the solution of the~\eqref{sde} will become almost sure exponentially stable.
\myendex

Using the random field approach described in Section~\ref{random:field} and considering~\eqref{heat:randomfields} driven by space-time white-noise, i.e.~$f=\delta_{0}$, one can show the following assertion regarding the moments of the mild random field solution $\{u(t,x) : t>0, x\in\mathbb{R}\}$. According to~\cite{FK,HUN}, for $k\geq 1$ and $C>0$, the moments grow exponentially, which means that
\begin{align*}
\mathbb{E}|u(t,x)|^{k}\sim\exp(Ck^{3}t),
\end{align*} 
 consequently
 \begin{align*}
 \tilde{\lambda}(k):=\limsup\limits_{t\to\infty}\frac{1}{t}\mathbb{E}\log|u(t,x)|^{k}\sim C k^{3}.
 \end{align*}
 Therefore, one observes
 \begin{align*}
 \tilde{\lambda}(1)<\frac{\tilde{\lambda}(2)}{2}< \ldots\frac{\tilde{\lambda}(k)}{k}<\ldots.
 \end{align*}
This phenomenon is called {\bf full intermittency}, which means that the random field solution develops high peaks concentrated on small sets for large time values. For further details on this topic, see~\cite{FK, HUN, Balan1, Balan2}.
Results regarding Lyapunov-exponents for parabolic SPDEs on {\em bounded domains} are available within the RDS approach using {\bf Oseledets' multiplicative ergodic theorem} for {\em compact operators} in~\cite{FlandoliS, Mohammed, CaraballoDuanLuSchmalfuss, lian, LuNS}. Beyond the Lyapunov spectrum, one is interested in further spectral properties for S(P)DEs, such as {\em dichotomy spectrum}. To our best knowledge this was dealt with only in the finite-dimensional setting~\cite{DSpectrum}.\\

If one is not interested in asymptotic stability as $t\ra +\I$, one may also study \textbf{metastability}, i.e., dynamics which looks stable on very long time scales. For example, consider the Allen-Cahn SPDE~\cite{BerglundAC}
\be
\label{eq:ACdynamics}
\txtd u =\left[ \partial_x^2 u + u-u^3\right] ~\txtd t + \sqrt{2\varepsilon}~\txtd W^\Id, u(0,x)=u_0(x),
\ee
where $u=u(t,x)$, $x\in\R/\Z=\mathbb{S}^1$, and we assume small noise $0<\varepsilon\ll1$. With the tools from Section~\ref{two:1}, one may show that~\eqref{eq:ACdynamics} has a global-in-time solution, e.g., a mild solution $u^\varepsilon_{u_0}(t)$, where the subscript reminds us of the dependence on the initial condition and the superscript of the parameter dependence. Furthermore, \eqref{eq:ACdynamics} is a stochastic {\bf gradient system} 
\beann
\txtd u &=& -\nabla_{\cH} V(u)~ \txtd t + \sqrt{2\varepsilon} ~\txtd W^\Id, \\
&&\text{where } V(\zeta):=\int_{\mathbb{S}^1}\frac12\|\nabla \zeta(x)\|^2_{\cH}+ \frac14(\zeta(x)^2-1)^2~\txtd x,
\eeann
where $V$ is called a {\bf potential} and where we use the space $\cH=L^2(\mathbb{S}^1)$. Let $P_t(h)(u_0):=\E[h(u^\varepsilon_{u_0}(t))]$, for a test function $h:\cH\ra \R$, denote the {\bf transition semigroup} associated to the solution process. In general, a measure $\mu$ is called an {\bf invariant measure} for $P_t$ if
\benn
P_t^*\mu=\mu\qquad \forall t\geq 0.
\eenn  
One may prove~\cite{DaPratoZabczyk} that quite a number of dissipative reaction-diffusion equations have an invariant measure using a tightness argument. For~\eqref{eq:ACdynamics}, one may write the invariant measure more explicitly as
\be
\label{eq:Gibbs}
\mu(\txtd \zeta)=\frac{1}{Z_0} \txte^{-V(\zeta)/\varepsilon}~\txtd \zeta,\quad Z_0:=\int_{\cH} \txte^{-V(\zeta)/\varepsilon}~\txtd\zeta,
\ee 
which is also known as the {\bf Gibbs measure} in this context; technically, the notation in~\eqref{eq:Gibbs} is to be understood more precisely using the Gaussian free field to get rid of the problem that there is no Lebesgue measure on $L^2(\cH)$. Although the invariant measure provides good insight regarding the bimodal stationary density of the gradient system~\eqref{eq:ACdynamics}, it does not provide direct information regarding finite-time dynamics for some fixed $T>0$ such as transition times between neighborhoods of deterministically stable steady states, e.g., between $u\equiv \pm1$. It is expected that the solution of~\eqref{eq:ACdynamics} approaches one of these equilibria and stays in a neighbourhood of these equilibria for a long time. Eventually, the system will transition from one of these neighborhoods to the other one. Such a behavior has been analyzed using large deviations in~\cite{BerglundGentzM}. The main ideas can be summarized as follows. We assume that $u^{\varepsilon}_{u_{0}}(t)$ converges to a deterministic quantity $\overline{u}_{u_{0}}(t)$ as $\varepsilon\to 0$, i.e., for any fixed $\delta>0$:
\begin{align*}
\lim\limits_{\varepsilon\to 0}\mathbb{P} \left(\sup\limits_{t\in[0,T]} \| u^{\varepsilon}_{u_{0}}(t)- \overline{u}_{u_{0}}(t) \|_{\cH} >\delta\right)=0.
\end{align*}
Furthermore, $u^{*}$ stands for an asymptotically stable equilibrium of the deterministic equation~\eqref{eq:ACdynamics}, $\cD_0:=\{u_{0}\in\cH :\lim\limits_{t\to\infty} \|\overline{u}_{u_{0}}(t) -u^{*} \|_{\cH}=0\}$ and
\begin{align*}
\tau^{\varepsilon}_{u_{0}} :=\inf\{ t>0 : u^{\varepsilon}_{u_{0}}(t) \notin \cD_0\}
\end{align*}
is the {\bf first-exit time}. As $\varepsilon\to 0$, the time required for this rare event to occur grows exponentially. Therefore, one investigates quantities like 
\benn
\lim\limits_{\varepsilon\to 0} \varepsilon \log \mathbb{E}\tau^{\varepsilon}_{u_{0}},\quad 
\lim\limits_{\varepsilon\to 0} \varepsilon \log \tau^{\varepsilon}_{u_{0}},\quad 
\lim\limits_{\varepsilon\to 0}u^{\varepsilon}_{u_{0}}(\tau^{\varepsilon}_{u_{0}}),
\eenn
which provide the expected exit-time, exit-time and exit location/shape.
For instance, one straightforward possibility to estimate the expected exit-time is 
\begin{align}\label{est:s}
\mathbb{E} (\tau^{\varepsilon}_{u_{0}}) \leq T \sum\limits_{k=0}^{\infty} \mathbb{P} (\tau^{\varepsilon}_{u_{0}}\geq kT ).
\end{align}
Due to the Markov property, one can show that for $k\in\mathbb{N}$:
\begin{align*}
\sup\limits_{u_{0}\in \cD_0 } \mathbb{P} (\tau^{\varepsilon}_{u_{0}} \geq kT )\leq \Big( \sup\limits_{u_{0}\in \cD_0 } \mathbb{P} (\tau^{\varepsilon}_{u_{0}} \geq T ) \Big)^{k}.
\end{align*}
Plugging this in~\eqref{est:s} leads to
\begin{align*}
\mathbb{E} (\tau^{\varepsilon}_{u_{0}})& \leq T \Big(1-\sup\limits_{u_{0}\in D} \mathbb{P}(\tau^{\varepsilon}_{u_{0}} \geq T) \Big)^{-1}\\
& \leq T \Big( \inf\limits_{u_{0}\in \cD_0 } \mathbb{P}(\tau^{\varepsilon}_{u_{0}}<T) \Big)^{-1}.
\end{align*}
This entails
\begin{align*}
\limsup\limits_{\varepsilon\to 0} \varepsilon \log\mathbb{E} (\tau^{\varepsilon}_{u_{0}}) \leq - \liminf\limits_{\varepsilon\to 0} \inf\limits_{u_{0}\in \cD_0 } \varepsilon \log \mathbb{P} (u^{\varepsilon}_{u_{0}} <T).
\end{align*}
In order to bound the right-hand side of the previous inequality, it is helpful to prove a {\bf large deviation principle (LDP)}~\cite{FreidlinWentzell,DaPratoZabczyk,FarisJona-Lasinio,Freidlin3,ChenalMillet,CerraiRoeckner,BerglundGentzM, HSS}. The key object in this case is the {\bf good rate function}
\benn
\cF_{[0,T]}(\gamma)=\frac12 \int_0^T \int_{\mathbb{S}^1} \left|\partial_t \gamma(t,x)\right|^2~\txtd x ~\txtd t,
\eenn
where $\gamma$ is a sufficiently regular path and we set $\cF_{[0,T]}(\gamma)=+\I$ if the last integral does not exist. One may prove $\cF_{[0,T]}$ is lower semi-continuous and has compact level sets. Furthermore, it helps us to estimate the probability of the solution process $(u_{u_0}^\varepsilon(t))_{t\in[0,T]}$ of~\eqref{eq:ACdynamics} as it satisfies an LDP
\beann
&&\liminf_{\varepsilon \ra 0} 2\varepsilon \ln \P\left((u_{u_0}^\varepsilon(t))_{t\in[0,T]}\in S_o\right)\geq -\inf_{\gamma\in S_o}\cF_{[0,T]}(\gamma),\\
&&\liminf_{\varepsilon \ra 0} 2\varepsilon \ln \P\left((u_{u_0}^\varepsilon(t))_{t\in[0,T]}\in S_c\right)\leq -\inf_{\gamma\in S_c}\cF_{[0,T]}(\gamma),
\eeann
where $S_o$ and $S_c$ denote the sets of all open and closed subsets of the space of continuous paths respectively. The idea is that since we are able via an LDP to control the probability that certain sample paths appear in the small noise regime, we can, e.g., analyze first-exit times~\cite{CerraiRoeckner,HSS} from a given sufficiently nice subset $\cD_0 \subset L^2(\cH)$. For example, the LDP implies for a subset $\Gamma$ of continuous paths that
\benn
\P\left((u_{u_0}^\varepsilon(t))_{t\in[0,T]}\in \Gamma\right)\approx \txte^{-\inf_\Gamma \cF_{[0,T]}/(2\varepsilon)}
\eenn  
to be understood asymptotically, and more precisely as logarithmic equivalence, as $\varepsilon\ra 0$. It is evident that the previous formulas indicate that {\bf metastability} occurs as the probabilities to exit from a deterministically stable steady state become exponentially small in the limit $\varepsilon\ra 0$. Using small noise limits and LDP estimates carries one quite far. The results range from sharp analytical exit-time asymptotics for reaction-diffusion systems~\cite{BerglundGentzM,Barret} to very broadly applicable numerical algorithms to capture metastability~\cite{KuehnSPDEcont,GrafkeSchaeferVandenEijnden}.  

\subsection{Invariant Manifolds}\label{im}

The aim of this section is to analyze sets that contain the trajectories of an SPDE that
converge to an equilibrium in forward / backward time or remain bounded for large time. For {\em deterministic} dynamical systems, these sets are called stable / unstable and center manifolds. Moreover, for initial conditions belonging to those sets,  the corresponding solution must also evolve within the set. This property is called {\bf invariance}. For deterministic systems, this can often be verified quite readily. For stochastic systems it is not a-priori clear, what a meaningful analogue of this concept is. We describe this fact within the RDS approach. We firstly recall that  $(\Omega,\mathcal{F},(\theta_{t})_{t\in\mathbb{R}},\mathbb{P})$ stands now for the {\em metric dynamical system}, constructed in Example~\ref{mds1}. The elements of this quadruple will be denoted by $\omega$ and $\theta_{t}\omega$ represents the shift introduced in~\eqref{shift}.\\ We illustrate the theory of invariant sets for stochastic evolution equations of the form~\eqref{pde:transf2} and indicate later on how this can be extended to more complicated SPDEs such as~\eqref{spde:fbm}. We focus on {\em center manifolds}, due to their numerous properties and applications, i.e.~center manifold theory allows a reduction to finite dimensional dynamics~\cite{DuDuan}. Therefore {\em stochastic center manifolds} have been intensively investigated in the literature, see~\cite{Arnold,DuDuan,WanngDuan, WaymireDuan, ChenRobertsDuan, ChekrounLiuWang, KN} and the references specified therein.

We impose the following restrictions on the coefficients of~\eqref{spde:transf2}. The drift   $f:\cH\to\cH$ is assumed to be a locally Lipschitz nonlinear term with $f(0)=f'(0)=0$. The spectrum of the linear operator $A$ is supposed to contain eigenvalues with zero and strictly negative real parts, i.e.~$\sigma(A)=\sigma^{\txtc}(A)\cup \sigma^{\txts}(A)$, where $\sigma^{\txtc}(A)=\{\lambda \in \sigma(A)\mbox{ : } \mbox{Re}(\lambda)=0\}$ and $\sigma^{\txts}(A)=\{\lambda\in \sigma(A) \mbox{ : } \mbox{Re}(\lambda)<0 \}$. The subspaces generated by the eigenvectors corresponding to these eigenvalues are denoted by $\cH^{\txtc}$ respectively $\cH^{\txts}$ and are referred to as {\em center} and {\em stable} subspace. These subspaces provide an invariant splitting of $\cH=\cH^{\txtc}\oplus \cH^{\txts}$.
We denote the restrictions of $A$ on $\cH^{\txtc}$ and $\cH^{\txts}$ by $A_{\txtc}:=A|_{\cH^{\txtc}}$ and $A_{\txts}:=A|_{\cH^{\txts}}$. Since $\cH$ is finite-dimensional we obtain that $S^{\txtc}(t):=\txte^{tA_{c}}$ is a {\em group} of linear operators on $\cH^{\txtc}$. Moreover, there exist projections $P^{\txtc}$ and $P^{\txts}$ such that $P^{\txtc} + P^{\txts} = \mbox{Id}_{\cH}$ and $A_{\txtc} =A|_{\cR(P^{\txtc})}$ and $A_{\txts}=A|_{\cR(P^{\txts})}$, where $\cR$ denotes the range of the corresponding projection. Additionally, we impose the following dichotomy condition on the semigroup. We assume that there exist two exponents $\tilde{\gamma}$ and $\tilde{\beta}$ with $-\tilde{\beta}<0\leq \tilde{\gamma}<\tilde{\beta}$ and constants $M_{c},M_{s}\geq 1$, such that
\begin{align}
&\|S^{\txtc}(t)  x\|_{\cH} \leq M_{c} \txte^{\tilde{\gamma} t} \|x\|_{\cH}, 
~~~\mbox{  for } t\leq 0 \mbox{ and } x\in \cH;\label{gamma:h}\\
& \|S^{\txts}(t) x\|_{\cH} \leq M_{s} \txte^{-\tilde{\beta} t} \|x\|_{\cH}, 
~~\mbox{for } t\geq 0 \mbox{ and } x\in \cH.\label{beta:h}
\end{align}
There are numerous operators that satisfy the spectral conditions imposed above. For instance let $\cH:=L^2(0,\pi)$ and set
$$A u:=\Delta u + u$$ with domain $\textnormal{Dom}(A)=H^2(0,\pi)\cap H^{1}_{0}(0,\pi)$.
Its spectrum is given by $\{1-n^{2}\mbox{ : } n\geq 1\}$ with corresponding eigenvectors $\{\sin(nx)\mbox{ : } n\geq 1\}$. The eigenvectors give us the center subspace $\cH^{c}=\mbox{span}\{\sin x\}$ and the stable one  $\cH^{s}=\mbox{span}\{\sin(nx)\mbox{ : } n\geq 2\}$.\\

We now investigate center manifolds for the SPDE~\eqref{spde:transf2},~under the above assumptions. Let $\varphi$ denote the RDS generated by~\eqref{spde:transf2}. 

\begin{defn}(\mbox{Random center manifold}) We call a set $\cM^{\txtc}(\omega)$ a {\bf random center manifold} if
	\begin{itemize}
		\item $\cM^{\txtc}(\omega)$ contains all trajectories $\varphi(t,\cdot,\cdot)$ which are bounded in forward and backward time;
		\item $\cM^{\txtc}(\omega)$ has a graph structure. This means that there exists a function $h^{\txtc }(\omega,\cdot):\cH^{\txtc}\to \cH^{\txts}$ with $h(\omega,0)=0$ such that
		\begin{align}\label{graph}
		\cM^{\txtc}(\omega) =\{ x + h^{\txtc}(\omega,x) : x\in \cH^{\txtc}\};
		\end{align}
		\item the {\bf tangency condition} $\txtD h^{\txtc}(\omega,0)=0$ holds true;
		\item $h^{\txtc}(\cdot,x):\Omega\to\cH^{\txts}$ is measurable for every $x\in\cH^{\txtc}$;
		\item $h^{\txtc}(\omega,\cdot):\cH^{\txtc}\to\cH^{\txts}$  is Lipschitz  / smooth;
		\item the following {\bf invariance} property holds true: if $x\in\cM^{\txtc}(\omega)$ then $\varphi(t,\omega,x)\in\cM^{\txtc}(\theta_{t}\omega)$ for $t\in\mathbb{R}^{+}$.
	\end{itemize}
If~\eqref{graph} is satisfied for $x\in \cH^{\txtc} \cap \cB(0,r(\omega))$, then $\cM^{\txtc}$ is called a {\em local} center manifold. Here $\cB(0,r(\omega))$ denotes a random neighbourhood of the origin.
\end{defn}

\begin{remark}
In the RDS approach the suitable concept for invariance of a random set~\cite{Arnold,DuanLuSchmalfuss} is that each orbit starting inside this random set, evolves and remains there omega-wise modulo the changes that occur due to the noise. These changes can be characterized by a suitable shift of the fiber of the noise. Another concept is constituted by the almost sure invariance of a deterministic set under stochastic influences, more precisely this means that each orbit starting inside this deterministic set remains there almost surely, see~\cite{WaymireDuan} and the references specified therein.
\end{remark}

We continue our deliberations regarding the existence of invariant sets for~\eqref{pde:transf2}. Since there are no stochastic differentials / integrals one can prove the existence of random center manifolds similar to the deterministic case. There one uses the Lyapunov-Perron method, which seeks that the trajectories of~\eqref{pde:transf2} that remain close to the center subspace under the dynamics. This can be equivalently formulated as a fixed-point problem in a suitable function space. More precisely, one introduces the continuous-time {\bf Lyapunov-Perron map / transform} for~\eqref{pde:transf2} is given by
 \begin{align}
 \label{lpeinfach}
 J(\omega,u,u(0))[t] & := S^{\txtc}(t) \txte^{\int\limits_{0}^{t} 
 	z(\theta_{\tau}{\omega}) ~\txtd \tau} P^{\txtc}u(0) \nonumber \\& 
 + \int\limits_{0}^{t} S^{\txtc}(t-r)  \txte^{\int\limits_{r}^{t} 
 	z(\theta_{\tau}{\omega}) ~\txtd \tau} P^{\txtc} 
 \overline{f}(\theta_{r}\omega, u(r))~\txtd r\nonumber\\
 & +\int\limits_{-\infty}^{t} S^{\txts}(t-r) \txte^{\int\limits_{r}^{t} 
 	z(\theta_{\tau}\omega) ~\txtd \tau} P^{\txts} \overline{f}(\theta_{r}\omega, 
 u(r))~\txtd r.
 \end{align}
 Further details regarding this operator can be found in~\cite{WanngDuan},~\cite[Section~6.2.2]{DuanWang},~\cite[Chapter~4]{ChekrounLiuWang} and the references specified therein. The next natural step is to show that~\eqref{lpeinfach} possesses a fixed-point in a certain function space. One possible choice turns out to be $BC^{\tilde{\eta},z}(\mathbb{R}^{-},\cH)$, see~\cite[p.~156]{DuanWang}. This space is defined as
 \begin{align*}
 BC^{\tilde{\eta},z}(\mathbb{R}^{-},\cH):=\left\{u:\mathbb{R}^{-}\to \cH, ~~
 \sup\limits_{t\leq 0}\txte^{-\tilde{\eta} t 
 	-\int\limits_{0}^{t}z(\theta_{\tau}\omega)~\txtd \tau } 
 \|u(t)\|_{\cH}<\infty\right\}
 \end{align*}
 and is endowed with the norm
 \begin{align}\label{bc}
 ||u||_{BC^{\tilde{\eta},z}} := \sup\limits_{t\leq 0}~\txte^{-\tilde{\eta} t 
 	-\int\limits_{0}^{t}z(\theta_{\tau}\omega)~\txtd \tau }\|u(t)\|_{\cH}.
 \end{align}
 Here $\tilde{\eta}$ is determined from~\eqref{gamma:h} and~\eqref{beta:h}, namely one has $-\tilde{\beta}<\tilde{\eta}<0$. Note that the previous expressions are well-defined since 
 \begin{align*}
 \lim\limits_{t\to \pm\infty}\frac{|z(\theta_t\omega)|}{|t|}=0, 
 \end{align*}
 according to~\cite[Lemma~2.1]{DuanLuSchmalfuss} and the references specified therein. Under a suitable smallness assumption on the Lipschitz constant of $f$ (gap condition) one can show that $J$ possesses a fixed-point $\Gamma(\cdot,\omega,u(0))$ for $u(0)\in \cH^{\txtc}$. Since such growth conditions on $f$ can be quite restrictive in applications, one usually introduces a cut-off function to truncate the nonlinearity outside a random ball around the origin. This fixed-point characterizes the random center manifold $\cM^{\txtc}(\omega)$ for~\eqref{pde:transf2}. More precisely, one can show that $\cM^{\txtc}(\omega)$ can be represented locally by the graph of a function $h^{\txtc}(\omega,\cdot)$, where $h^{\txtc}(\omega,u(0))=P^{s}\Gamma(0,\omega,u(0))$ for $ u(0) \in \cH^{\txtc} \cap \cB(0,r(\omega))$, i.e.
 \begin{align}\label{h:strat}
 h^{\txtc}(\omega,u(0)) = \int\limits_{-\infty}^{0} S^{\txts}(-\tau) \txte^{\int\limits_{\tau}^{0}z(\theta_{r}\omega)~\txtd r} P^{\txts} \overline{f} (\theta_{\tau}\omega,\Gamma(\tau,\omega,u(0)))~\txtd \tau.
 \end{align}
 Here $\cB(0,r(\omega))$ denotes a random neighbourhood of the origin, i.e.~the radius $r(\omega)$ depends on the intensity/magnitude of the noise.
 \begin{ex}
  Let $a>0$, $\sigma>0$ and $W$ stand for a one-dimensional Brownian motion. For the SPDE
 	\begin{align}\label{ex1}
 	\begin{cases}
 	\txtd u = (\Delta u + u - a u^{3}) ~ dt + \sigma u \circ \txtd W(t) \\
 	u(0,t)=u(\pi,t)=0, ~\mbox{for } t \geq 0\\
 	u(x,0)=u_{0}(x),~~~~~~~\mbox{for  } x\in[0,\pi],
 	\end{cases}
 	\end{align}
 	the transformation into a PDE with random coefficients leads to
 	 \begin{align}\label{random:pde}
 	\frac{\partial u}{\partial t} = \Delta u + u + \sigma z(\theta_{t}\omega) u - a \txte^{2\sigma z (\theta_{t}\omega)}u^{3},
 	\end{align}
 	as discussed in Example~\ref{transf2}.
 		Regarding the discussion above, one can infer that~\eqref{ex1} has a local center manifold
 	$$\cM^{\txtc}(\omega) =\{B\sin x + h^{\txtc}(\omega, B \sin x)\}=\Bigg\{ B\sin x +\sum\limits_{n=2}^{\infty} c_{n}(\omega,B)\sin (nx) \Bigg\}.$$
 	In this case, it is also possible to derive suitable approximation results for $h^{\txtc}$, namely one can show that the coefficients satisfy $c_{n}(\omega,B)=\mathcal{O}(B^{3})$ as $B\to 0$.
 	Plugging this in~\eqref{random:pde} gives us the {\em nonautonomous} equation on the {\em center manifold}
 	\begin{align*}
 	\frac{\txtd}{\txtd t} (B \sin x)= B\sigma z (\theta_{t}\omega) \sin x - \frac{3}{4} a B^{3}\sin x \txte^{2\sigma z(\theta_{t}\omega)} + \cO(B^{5}),
 	\end{align*}
 	consequently
 	\begin{align*}
 	\dot{B} = \sigma z(\theta_{t}\omega) B - \frac{3}{4} a B^3 \txte^{2\sigma z (\theta_{t}\omega)}+  \cO(B^{5}).
 	\end{align*}
 	Since $-u$ is also a solution for~\eqref{random:pde} we have that $c_{n}(\omega,B)=0$ for $n$ even. Therefore one has the following approximation of $h$ as
 	\begin{align*}
 	h^{\txtc}(\omega, B\sin x) = c_{3}(\omega,B)\sin 3x + \cO(B^{5}).
 	\end{align*}
 	Consequently, all the expressions arising in the Lyapunov-Perron method can be explicitly computed in this situation.
 \myendex 

The next step is to establish the existence of invariant sets for SPDEs {\em without} reducing them into PDEs with random coefficients. Using the rough path techniques presented in Section~\ref{rough:paths}, one can immediately infer that the solution operator of~\eqref{spde:fbm} generates a random dynamical system, see~\cite{HesseNeamtu2}. Since the stochastic convolution / integral~\eqref{integral} was constructed in a pathwise way and not {\em almost surely}, there are no exceptional sets that can occur. Therefore, one can solve the SPDE~\eqref{spde:fbm} for every random input $\omega$ which is $\alpha$-H\"older continuous, for $\alpha\in(1/3,1/2]$. This case includes Brownian motion and fractional Brownian motion with Hurst index $H\in(1/3,1/2)$, as illustrated in Section~\ref{rough:paths}. Consequently, the results obtained using the transformation of an SPDE with linear multiplicative noise into a random PDE discussed above, will be recovered in a more general setting.\\
Under the assumptions on the linear operator $A$ and drift term $f$ specified above we can now investigate random sets for the SPDE with multiplicative noise~\eqref{spde:fbm}. Using the notations in Section~\ref{rough:paths} we additionally assume that $g(0)=0$. Regarding the deliberations above and~\eqref{lpeinfach}, we infer that the {\bf Lyapunov-Perron map/transform} is given by
\begin{align}
&J(\omega,u,u(0))[t]: = S^{\txtc} (s) P^{c} u(0) + \int\limits_{0}^{t}S^{\txtc}(t-r) P^{\txtc} f(u(r))~\txtd r \nonumber \\
& +\int\limits_{0}^{t} S^{\txtc} (t-r) P^{\txtc}g(u(r))~\txtd \omega(r) + \int\limits_{-\infty}^{t} S^{\txtc} (t-r) P^{\txts} f(u(r))~\txtd r \nonumber\\
& +\int\limits_{-\infty}^{t} S^{\txts}(t-r)P^{\txts}g(u(r))~\txtd\omega(r),\label{lp:rp}
\end{align}
where for notational consistency we used the identification $B^{H}(t,\omega)=\omega(t)$. Unlike~\eqref{lpeinfach}, the Lyapunov-Perron transform $J$ contains stochastic integrals and it is not clear in which function space one should formulate the fixed-point problem. However, modifying the Lyapunov-Perron method one obtains the following result, see~\cite{KN}.
\begin{thm}
	There exists a local center manifold $\cM^{\txtc}(\omega)$ for~\eqref{spde:fbm} such that 
	\begin{align*}
	\cM^{\txtc}(\omega) =\{u(0) + h^{\txtc}(\omega,u(0)) : u(0)\in \cH^{\txtc}\cap \cB(0,r(\omega))\},
	\end{align*}
	where 
	\begin{align*}
	h^{\txtc}(\omega,u(0)) &= \int\limits_{-\infty}^{0} S^{\txts}(-r)P^{\txts} f(\Gamma(r,\omega,u(0)))~\txtd r \\
	&+ \int\limits_{-\infty}^{0} S^{\txts}(-r)P^{\txts} g(\Gamma(r,\omega,u(0)))~\txtd \omega(r)
	\end{align*}
	and $\Gamma$ is the fixed-point of $J$.
\end{thm}

We are going to outline the ideas of the proof. Since the estimates of the stochastic integrals appearing in the definition of $J$ contain certain H\"older norms of the random input $\omega$, which are uniform over the unit interval, it is meaningful to:
\begin{itemize}
	\item [1)] discretize the Lyapunov-Perron map~\eqref{lp:rp}, i.e.~consider its solution at {\em discrete} time-points and obtain a sequence of solutions over $[0,1]$;
	\item [2)] apply the Lyapunov-Perron method in a space of {\em sequences}, i.e.~formulate the fixed-point problem for the discrete version of $J$ in a space of sequences (instead of $BC^{\eta,z}$ as above).
\end{itemize}
The proof essentially combines rough path techniques with the Lyapunov-Perron method for {\em discrete-time} dynamical systems. For further details on this topic see~\cite{KN}.
\begin{remark}
To obtain stable / unstable manifolds, one imposes different spectral assumptions on the linear part and modifies the definition of $J$ accordingly, see~\cite{DLS,DuanLuSchmalfuss,GarridoLuSchmalfuss, LS}.
\end{remark}
We conclude this section pointing out once more that the time-evolution of a manifold in the RDS framework is described using an appropriate shift with respect to the noise. However, there are further theories that provide bounds on the probabilities that these manifolds evolve in time, see for instance~\cite{BloemkerHairer,Bloemker,BloemkerWang}.\\
Similar to center manifold theory that often allows for a local reduction to finite dimensional dynamics for SPDEs, there are several other approaches that provide approximations for the solutions of SPDEs by finite-dimensional SDEs as in~\cite{BloemkerHairer, BloemkerMohammed}. To briefly illustrate such results, we consider similar to Example~\ref{ex1}, for small $\varepsilon>0$ and arbitrary $\nu,\sigma>0$ the SPDE
\begin{align}\label{spde:e}
\txtd u  = [(\Delta u + u) +\nu\varepsilon^2 - u^3 ]~\txtd t+ \sigma\varepsilon~\txtd W(t),
\end{align}
with Dirichlet boundary conditions on $[0,\pi]$. The term $\nu\varepsilon^2$ is a small linear perturbation. As seen in Example~\ref{ex1}, $\cH^{\txtc}=\mbox{span}\{\sin x\}$ and $\cH^{\txts}=\mbox{span}\{\sin (nx) : n\geq 2\}$. Assuming that the additive noise is acting for instance only on $\sin(2x)$ and rescaling time as $T:=\varepsilon^2 t$, one can show that the rescaled solution of the SPDE~\eqref{spde:e}
\begin{align*}
u(t):=\varepsilon v (\varepsilon^2 t),
\end{align*}
can be approximated as
\begin{align*}
v(T)\approx B(T)\sin x + \frac{\sigma}{3} \tilde{Z}(T) \sin (2x) + \mathcal{O}(\varepsilon^{1-}).
\end{align*}
Here $\mathcal{O}(\varepsilon^{1-})$ collects the higher-order terms, $B$ is the solution of a certain {\bf amplitude equation} (or {\bf modulation equation}) given by
\begin{align*}
\partial_T B = \Big( \nu-\frac{\sigma^2}{4} \Big) B - \frac{3}{4} B^{3}
\end{align*}
and $\tilde{Z}$ is the Ornstein-Uhlenbeck process on the new time-scale given by
\begin{align*}
\tilde{Z}(T)= \varepsilon^{-1}\int\limits_{0}^{T} \txte^{-\lambda_2\varepsilon^2(T-\tau)}~\txtd \tilde{W}(\tau),
\end{align*}
where $\lambda_{2}$ is the eigenvalue corresponding to $\sin(2x)$ and $\tilde{W}(T):=\varepsilon W({\varepsilon^{-2}T})$, is a rescaled Brownian motion. Naturally, if the additive noise acts on the other modes of $\cH^{\txts}$, i.e.~$W(t)=\sum\limits_{k=2}^{N}\sigma_{k}w_{k}(t)e_{k}$ for a fixed $N>2$ or $W(t)=\sum\limits_{k=2}^{\infty}\sigma_{k}w_{k}(t)e_{k}$ together with a decay condition on the coefficients of this series, one can derive analogous amplitude equations and approximation results. For the proof of such statements and further applications, see~\cite{Bloemker,BloemkerHairer,BloemkerMohammed} and~\cite[Section~6.6]{WanngDuan}.

\subsection{Random Attractors}\label{ra}

As already emphasized, RDS are generalizations of nonautonomous dynamical systems. The next concept is a generalization of a {\em deterministic attractor} for nonautonomous systems, as considered in~\cite{BabinVishik, Temam, SellYou, KloedenRasmussen, KPR}. Intuitively, an attractor is a compact set of the phase space towards which the dynamical system evolves after a certain amount of time. Before dealing with attractors for SPDEs, we give a concise example which illustrates, why is it meaningful for nonautonomous systems to take a {\em pullback limit} with respect to {\em time} instead of a {\em forward} one. 
\begin{ex}
Let $a, b>0$ and consider as in~\cite{Ghil} the simple nonautonomous ODE
	\begin{align*}
	\begin{cases}
	u' = -a u + b t\\
	u(0)=u_{0}.
	\end{cases}
	\end{align*}
	Its solution $u(t)= \txte^{-at}u_0 +\frac{b t}{a} -\frac{b}{a^{2}} $ obviously blows up for $t\to\infty$. However, if one considers the solution operator at time $t$ with $u(s)=u_0$ for $t\geq s$, i.e.~$\overline{\varphi}(t,s,u_0)$, and takes the pullback limit $s\to-\infty$ for a fixed $t>0$, a straightforward computation entails that
	\begin{align}\label{attr:1}
	|\overline{\varphi}(t,s,u_{0}) - \overline{\mathcal{A}}(t)|\to 0, \mbox{ as } s\to-\infty,
	\end{align}
where $\overline{\cA}(t):= \frac{b t}{a} -\frac{b}{a^{2}}$.	
\myendex 
 For an SPDE, an attractor will be a {\em random, time-dependent} set. We present the mathematical formalism of this concept and explain what {\em pullback} in this framework means.~Again, we assume that the solution operator of an SPDE generates a random dynamical system $\varphi$ on $\cH$ and present concrete examples later on.~We recall that throughout this section $(\Omega,\mathcal{F},(\theta_{t})_{t\in\mathbb{R}},\mathbb{P})$ denotes the metric dynamical system constructed in Example~\ref{mds1}.

\begin{defn} 
	A random bounded set $\{\cS(\omega)\}_{\omega\in\Omega}$ of $\cH$ is called {\bf tempered} with respect to $(\theta_{t})_{t\in\mathbb{R}}$ if for all $\omega\in\Omega$ it holds
	$$ \lim\limits_{t\to\infty}\txte^{-\overline{\beta} t}\sup\limits_{x\in \cS(\theta_{-t}\omega)}\|x\|=0,\mbox{  for all } \overline{\beta}>0.$$
\end{defn}

In the sequel, $\mathcal{D}$ denotes the collection of tempered random sets in $\cH$.

\begin{defn}(\mbox{Random absorbing set})\label{absorbing}
		A set $\{\cS(\omega)\}_{\omega\in\Omega}\in\mathcal{D}$ is called {\bf random absorbing set} for $\varphi$ if for every $\cR=\{\cR(\omega)\}_{\omega\in\Omega}\in\mathcal{D}$ and  $\omega\in\Omega$, there exists a random  time $t_{\cR}(\omega)>0$ such that
	$$ \varphi(t,\theta_{-t}\omega,\cR(\theta_{-t}\omega))\subseteq \cS(\omega),\mbox{ for all } t\geq t_{\cR}(\omega).$$
\end{defn}

\begin{defn}(\mbox{Random pullback attractor})
		A random set $\{\mathcal{A}(\omega)\}_{\omega\in\Omega}\in\mathcal{D}$ is called a 
	$\mathcal{D}$-{\bf random (pullback) attractor} for $\varphi$ if the following properties are satisfied:
	\begin{itemize}
		\item [a)] $\mathcal{A}(\omega)$ is compact for every $\omega\in\Omega$;
		\item [b)] $\{\mathcal{A}(\omega)\}_{\omega\in\Omega}$ is positive invariant, i.e.
		$$\varphi(t,\omega,\mathcal{A}(\omega))=\mathcal{A}(\theta_{t}\omega)\mbox{ for all } t\geq 0;$$
		\item [c)] $\{\mathcal{A}(\omega)\}_{\omega\in\Omega}$ pullback attracts every set in $\mathcal{D}$, more precisely, for every $\cR=\{\cR(\omega)\}_{\omega\in\Omega}\in\mathcal{D}$,
		\begin{align}\label{attr:a}
		\lim\limits_{t\to\infty}\mbox{dist}_{\cH}(\varphi(t,\theta_{-t}\omega,\cR(\theta_{-t}\omega)),\mathcal{A}(\omega))=0,	
		\end{align}
		where $\mbox{dist}_{\cH}$ is the {\bf Hausdorff semi-distance}. This is given by
		$\mbox{dist}_{\cH}(Y,Z)=\sup\limits_{y\in Y}\inf\limits_{z\in Z}\|y-z\|$, for any subsets $Y\subseteq \cH$ and $Z\subseteq \cH$.
	\end{itemize}
	If the set $\cA(\omega)$ satisfies all the properties above and consists only of a singleton, then this is called {\bf random fixed-point}.
	
\end{defn}
From this abstract definition, we see that the concept pullback means in this setting, that one shifts into the past the fiber of the noise, compare~\eqref{attr:1} and~\eqref{attr:a}. For a better comprehension we present the following example.
\begin{ex} Let $\overline{\mu}>0$ and consider the one-dimensional SDE
\begin{align}\label{lin:sde}
\txtd u= -\overline{\mu} u~\txtd t + \txtd W(t),
\end{align} 
Obviously, for $t\geq t_{0}$, its solution is given by
\begin{align*}
u(t) =u(t_{0}) \txte^{-\overline{\mu}(t-t_{0})} + \txte^{-\overline{\mu} t}\int\limits_{t_{0}}^{t} \txte^{\overline{\mu}s}~\txtd W(s).
\end{align*}
Note that the forward limit $t\to \infty$ does not exist. However, for a fixed $t\in\mathbb{R}$ taking the pullback 
limit $t_{0}\to-\infty$  yields the stationary solution of this SDE
\begin{align*}
\lim\limits_{t_{0}\to-\infty} u(t) =  \int\limits_{-\infty}^{t} \txte^{-\overline{\mu}(t-s)}~\txtd W(s)=\int\limits_{-\infty}^{0} \txte^{\overline{\mu}s}~ \txtd \theta_{t}W(s),
\end{align*}
which is the one-dimensional Ornstein-Uhlenbeck process $(t,\omega)\mapsto z(\theta_{t}\omega)$
$$z(\theta_{t}\omega) =\int\limits_{-\infty}^{0} \txte^{\overline{\mu}s} ~\txtd\theta_{t}\omega(s).$$
 One can show that $z$ is the random fixed-point of~\eqref{lin:sde}, i.e.~in this case the attractor $\mathcal{A}(\omega):=\{z(\omega)\}$ is a singleton.~The same statement holds true in the infinite-dimensional setting, recall~\eqref{lin:spde} and~\eqref{ou}.
\myendex 

The existence of random attractors can be shown by the following criterion, see Theorem 2.1 in~\cite{schmalfuss} and~\cite{CrauelFlandoli,crauel}.
\begin{thm}\em{(Criterion for existence of a pullback attractor)}\label{ex:attractor}
		Let $\varphi$ be a continuous random dynamical system on $\cH$ over $(\Omega,\mathcal{F},\mathbb{P},(\theta_{t})_{t\in\mathbb{R}})$. 
	Suppose that $\{\cK(\omega)\}_{\omega\in\Omega}$ is a compact random absorbing set for $\varphi$ in $\mathcal{D}$. 
	Then $\varphi$ has a unique $\mathcal{D}$-random attractor $\{\mathcal{A}(\omega)\}_{\omega\in\Omega}$ which is given by
	$$\mathcal{A}(\omega)=\bigcap\limits_{\tau\geq 0}\overline{\bigcup\limits_{t\geq\tau}\varphi(t,\theta_{-t}\omega,\cK(\theta_{-t}\omega))}. $$
\end{thm}

This theorem indicates that one has to verify the following two aspects in order to prove the existence of a random pullback attractor.
\begin{itemize}
	\item [1)] The existence of an absorbing set. This usually follows by suitable a-priori estimates on the solutions. More precisely, the following condition is convenient to show the existence of an absorbing set. If for every $x\in \cR(\theta_{-t}\omega)$, $\cR\in\mathcal{D}$ and $\omega\in\Omega$ it holds
	\begin{align}\label{criterion}
	\limsup\limits_{t\to\infty} \|\varphi(t,\theta_{-t}\omega,x)\|_{\cH}\leq \rho(\omega),
	\end{align}
	where $\rho(\omega)>0$ is a {\em tempered random variable}, then the ball centered in $0$ with radius $\rho(\omega)+\delta$, 
	i.e.~$\cS(\omega):=\mathcal{B}(0,\rho(\omega)+\delta)$, for some constant $\delta>0,$ is a random absorbing set.
	\item [2)]  The compactness of the absorbing set. This follows in general by deriving estimates of the solutions in function spaces that are compactly embedded in $\cH$.
\end{itemize}
For further details and applications see~\cite{schmalfuss, CrauelFlandoli}. 
We now discuss how these abstract results can be applied to concrete SPDEs, more precisely to SPDEs which are equivalent to~\eqref{pde:transf1}. In order to obtain the necessary compactness, we additionally assume that the semigroup $(S(t))_{t\geq 0}$ is analytic and {\em compact}. Furthermore, we assume that there exist constants $\widetilde{M}\geq 1$ and $\overline{\mu}>0$ such that $\|S(t)\|_{\cH}\leq \widetilde{M}\txte^{-\overline{\mu}t}$ for $t>0$. Therefore, the stationary Ornstein-Uhlenbeck process
$$Z(\theta_{t}\omega) =\int\limits_{-\infty}^{t} S(t-s)~\txtd \omega(s) $$
is well-defined. An example in this sense is given by the equation
\begin{align*}
\txtd u =[ \Delta u - \overline{\mu} u]~\txtd t + f(u+Z(\theta_{t}\omega))~\txtd t
\end{align*}
on $\cH:=L^{2}(D)$, where $D\subset\mathbb{R}^d$ is an open bounded domain.
The most simple situation is to assume the nonlinear term $f:\cH\to\cH$ is globally Lipschitz continuous with Lipschitz constant $L_{f}>0$ and of bounded growth with constant $0<l_{f}<L_{f}$. In this case, letting $\varphi(t,\omega,u_{0})$ denote the corresponding solution operator of~\eqref{pde:transf1}, the Gronwall inequality immediately entails the a-priori bound
\begin{align*}
\|\varphi(t,\omega,u_{0})\|_{\cH} \leq \widetilde{M} \txte^{(l_{f}\widetilde{M} - \overline{\mu}) t }\|u_{0}\|_{\cH} + M \int\limits_{0}^{t}C(\theta_{s}\omega) e ^{(l_{f}\widetilde{M}-\overline{\mu})(t-s) }~\txtd s,
\end{align*}
where $C(\omega):=l_{f}\|Z(\omega)\|_{\cH} +c$ and $c>0$ is a constant. This terms occurs due to the structure and growth boundedness of $f$. Replacing $\omega$ by $\theta_{-t}\omega$ in the previous inequality and further assuming that $\overline{\mu}-l_{f}\widetilde{M}>0$, leads to
\begin{align*}
\|\varphi(t,\theta_{-t}\omega,u_{0})\|_{\cH} &\leq \widetilde{M} \txte^{(l\widetilde{M}-\overline{\mu}) t} \|u_{0}\|_{\cH} + \int\limits_{-t}^{0} C(\theta_{\tau}\omega) \txte^{(\overline{\mu} -l_{f}\widetilde{M})\tau }~\txtd\tau\\
& \leq  \widetilde{M} \txte^{(l\widetilde{M}-\overline{\mu}) t} \|u_{0}\|_{\cH} + \int\limits_{-\infty}^{0} C(\theta_{\tau}\omega)\txte^{(\overline{\mu} -l_{f}\widetilde{M})\tau}~\txtd\tau.
\end{align*}
This bound immediately entails the existence of an absorbing set if  $\cS(\omega):=\mathcal{B}(0,\rho(\omega)+\delta)$, where $\delta>0$ and the tempered random variable $\rho$ is given by
\begin{align*}
\rho(\omega):= \int\limits_{-\infty}^{0} C(\theta_{\tau}\omega)\txte^{(\overline{\mu} -l_{f}\widetilde{M})\tau}~\txtd\tau.
\end{align*}
Keeping this in mind, a natural candidate for a {\em compact} absorbing set for $\varphi$ would be for instance $\cK(\omega):=\overline{\varphi(t^*,\theta_{-t^{*}}\omega, \cS(\theta_{-t^{*}}\omega))}$, where $t^{*}>0$ is fixed and the closure is taken with respect to the topology in $\cH$.
This can be achieved showing that $$\|\varphi(t^*,\theta_{-t^{*}}\omega, \cS(\theta_{-t^{*}}\omega)) \|_{\cH_{\gamma}}<\infty,$$ where $\gamma\in(0,1)$ and $\cH_{\gamma}=\textnormal{Dom}((-A)^{\gamma})$ are the domains of the fractional powers of the operator, as discussed in Section~\ref{rough:paths}. Due to the fact that $(S(t))_{t\geq 0}$ is a compact $C_{0}$-semigroup, the embedding $\cH_{\gamma}\hookrightarrow\cH$ is compact which proves the compactness of the random set $\cK(\omega)$. Due to Theorem~\ref{ex:attractor} we infer that~\eqref{pde:transf1} has a random pullback attractor.\\
Similar arguments can be employed if the nonlinear term satisfies suitable dissipativity conditions. More precisely, we consider the following equation
\begin{align}\label{attr:bd}
\txtd u = [\Delta u - \overline{\mu} u]~\txtd t + f(x,u +Z(\theta_{t}\omega))~\txtd t, ~~t>0\mbox{ and } x\in D.
\end{align}
As commonly met in the theory of reaction-diffusion equations, see~\cite[Section~5.1]{SellYou}, we impose the following restrictions on the reaction terms. For $x\in D$, $s\in\mathbb{R}$ and an integer $q\geq 2$ we assume that there exist positive constants $\alpha_{1},\alpha_{2},\alpha_{3}$ and $c_{1},c_{2},c_{3}$ such that
\begin{align*}
&f(x,s)s\leq -\alpha_{1}|s|^{q} + c_{1}\\
&|f(x,s)|\leq\alpha_{2}|s|^{q-1} + c_{2}\\
& \Big| \frac{\partial f}{\partial s}(x,s) \Big|\leq \alpha_{3} |s|^{q-2}+c_{3}.
\end{align*}
This means that $f$ can be a polynomial of odd degree with a negative leading order coefficient. In this case one can derive a-priori estimates of the solution of~\eqref{attr:bd} in $L^{2}(D)$ to obtain the existence of an absorbing set, and in $W^{1,2}(D)$ for the compactness argument.\\
Major technical difficulties occur when compact embeddings are no longer available. This is for instance the case, when one considers~\eqref{attr:bd} on unbounded domains, i.e.~on $\mathbb{R}^d$, see~\cite{Bates}. In this case, one replaces the constants $c_{1},c_{2},c_3$ in the previous assumption with $L^{2}(\mathbb{R}^d)$-integrable functions and modifies~\eqref{attr:bd} as
\begin{align}\label{attr:ubd}
\txtd u = [\Delta u - \overline{\mu} u]~\txtd t + [f(x,u +Z(\theta_{t}\omega)) + h(x) ]~\txtd t, 
\end{align}
for $t>0$ and $x\in\mathbb{R}^d$, where $h\in L^{2}(\mathbb{R}^d)$. Here, one needs again estimates in $L^{2}(\mathbb{R}^{d})$ and $W^{1,2}(\mathbb{R}^d)$. The most technical argument is to show by means of a cut-off technique that the tails of the solution of~\eqref{attr:ubd} will become uniformly small for large enough time. Afterwards one can use the compact embedding $W^{1,2}(D)\hookrightarrow L^{2}(D)$ for bounded domains $D\subset\mathbb{R}^d$ and the uniform decay for sufficiently large time outside of $D$. The function $h$ is required in order to show that the attractor for~\eqref{attr:ubd} is set-valued. Without $h$, the attractor reduces to a singleton.\\
Another situation where compact embeddings cannot be directly employed occurs in the context of {\bf partly dissipative systems}, i.e.~coupled SPDEs with SDEs. An example in this sense is given by the following system. We assume for simplicity that $x\in [0,1]$ and consider
\begin{align*}
 &\txtd u_1 = \Delta u_{1}~\txtd t + [f_1(x,u_1+Z_{1}(\theta_{t}\omega, u_2+Z_{2}(\theta_{t}\omega))) +f_{2}(x,u_{1}+Z_{1}(\theta_{t}\omega)) ]~\txtd t\\
& \txtd u_{2} =-\nu u_{2} ~\txtd t+ f_{3}(x,u_{1}+Z_{1}(\theta_{t}\omega))~\txtd t,
\end{align*}
where $\nu>0$, $Z_{1},Z_{2}$ are two Ornstein-Uhlenbeck processes and $f_{1},f_{2},f_{3}$ satisfy suitable dissipativity assumptions. The technical difficulty consists in the missing regularizing effect of the Laplacian in the second component. To overcome this, one splits the second component $u_2$ into a regular part (i.e.~which belongs to $W^{1,2}(0,1)$) and a remaining part which tends asymptotically to zero. This is sufficient to obtain the necessary compactness results. Further details on this topic can be looked up in~\cite{KuehnNPein}.\\ 
Before pointing out some concluding remarks we would like to emphasize that the structure of random attractors can be totally different from the deterministic case. For instance, consider as in~\cite{FlandoliGessScheutzow}, the porous-media equation
\begin{align}\label{pom}
\txtd u = [\Delta |u|^{m-1} u + u]~\txtd t +\txtd W^{Q}(t),
\end{align}
with zero Dirichlet boundary conditions on a bounded domain $D\subseteq\mathbb{R}^d$, where $d\leq 4$ and $m>1$. The noise $W^{Q}$ is a trace-class Brownian motion on $W^{-1,2}(D)=(W^{1,2}_{0}(D))^{*}$. It is known that the attractor of the {\em deterministic} porous-media equation has infinite fractal dimension. However, the {\em random} attractor of~\eqref{pom} reduces to a fixed-point, see~\cite[Section~3]{FlandoliGessScheutzow}. This is an example of {\em synchronization by noise} as investigated in~\cite{FlandoliGessScheutzow}.
\begin{remark}
	\begin{itemize}
		\item [1)] Results on the existence of random attractors for SPDEs without directly reducing the SPDE in a random PDE are available in~\cite{BessaihGarridoSchmalfuss}. Here one deals with semilinear {\em delay equations} under certain smoothness assumptions on the coefficients. These allow an {\em integration by parts formula} which is employed in order to get rid of the stochastic integrals and do a pathwise analysis of the equation. Apart from this, all results regarding random pullback attractors for SPDEs with {\em additive} or {\em linear multiplicative} noise employ transformations in a PDE with random coefficients as described above~\cite{attrNS, Bates,CrauelFlandoli, FlandoliGessScheutzow, Gess1,Gess2,GessLiuRoeckner, Ghil, schmalfuss, Wang}.~To our best knowledge, for nonlinear multiplicative noise, and SPDEs such as~\eqref{spde:fbm}, there are no general existence results concerning random attractors. For multiplicative fractional noise, in the more regular case, i.e.~$H\in(1/2,1)$, results on this topic are contained in~\cite{Gao}.
		\item [2)] For assertions regarding the dimension of random attractors, see~\cite{Debussche, LangaRobinson} and more recently~\cite{CaSo}. 
		\item [3)] There are several further concepts of random attractors, for instance {\em weak}~\cite{Scheutzow:a,FlandoliGessScheutzow}, {\em exponential}~\cite{CaSo} or {\em mean-square} random attractors~\cite{KloedenLorenz}.
	\end{itemize}
\end{remark}

\subsection{Further Dynamical Aspects}

In this section we just point out several dynamical phenomena for SPDEs we have not discussed above due to practical space limitations within this edited volume. We also provide some references for the reader to get started in the background literature on these topics.~We emphasize that some of the following aspects have been investigated especially for additive and multiplicative white-in-time noise while the influence of other type of time-correlated noise is highly challenging and remains to be explored in far more detail. 
Such dynamical phenomena include among others: 

\begin{itemize}
	\item traveling waves and their stability~\cite{Kuehn1,MuellerMQ,HH};
	\item pattern formation~\cite{GarciaOjalvoSancho,GKuehn,Huttetal};
	\item stochastic bifurcations~\cite{Arnold,BHP,KuehnRomano1};
	\item early-warning signs for bifurcations in SPDEs~\cite{GKuehn,KuehnRomano1,KuehnFKPP};
	\item sample paths estimates for fast-slow SPDEs~\cite{BerglundGentzM,GnannKuehnPein};
	\item averaging results for fast-slow SPDEs~\cite{WangRoberts,DuanWang};
	\item stochastic homogenization~\cite{Cerrai,FW,BessaihEM}.
	\item finite-time blow-up~\cite{MuellerBlowUp,LvDuan}.
\end{itemize}

Naturally, there exist further probabilistic methods  that provide dynamical insights, which are based on Kolmogorov/Fokker-Planck equations associated to SPDEs~\cite{Bogachev,BarbuRoeckner, DaPratoZabczyk}.
For higher-order-in-time SPDEs, such as wave equations / dispersive equations there are available results regarding invariant manifolds~\cite{LS} and random attractors~\cite{Wang,CaSo}.\\
In summary, a lot of progress towards solution theory of SPDEs has been made over the last several decades. The development of dynamical system results for SPDEs is currently actively growing and promises to be an active area of research for many decades to come. For example, our understanding of the structure, formation, stability, and interaction of patterns for SPDEs is still far away from the level of results available in the PDE setting. A similar remark applies to bifurcation problems and the role played by singular SPDEs in this context.\medskip

\textbf{Acknowledgements}. CK and AN have been supported by a DFG grant in the D-A-CH framework (KU 3333/2-1). CK acknowledges support by a Lichtenberg Professorship.~CK and AN would like to thank Prof.~Wilfried Grecksch and Prof.~Hannelore Lisei for the kind invitation to contribute to the book {\em Finite and Infinite Dimensional
	Stochastic Equations with Applications in Physics}.

\end{document}